\documentclass[12pt]{article}

\usepackage{graphicx,psfrag,subfigure,amsmath}
\usepackage{url}
\input epsf
\usepackage{latexsym}
\usepackage{amsfonts, amssymb}
\usepackage{mathrsfs,color}
\newtheorem{thm}{Theorem}[section]
\newtheorem{defi}{Definition}[section]
\newtheorem{claim}{Claim}[section]

\newtheorem{lem}{Lemma}[section]

\newcommand\qed{\hfill $\Box$ \smallskip \smallskip}

\newenvironment{wst}
{\setlength{\leftmargini}{1.5\parindent}
 \begin{itemize}
 \setlength{\itemsep}{-1.1mm}}
{\end{itemize}}

\begin{document}




\title{\bf Erd\H{o}s-Gy\'{a}rf\'{a}s Conjecture for $P_{10}$-free Graphs\thanks{Supported by  NSFC grants 11771172, 11871239 and 11971196.}}

\author{\small Zhiquan Hu\thanks{School of Mathematics and Statistics,   Central China Normal
University,  Wuhan, China, hu\_zhiq@aliyun.com},
\, Changlong Shen\thanks{School of Mathematics and Statistics,   Central China Normal
University,  Wuhan, China, clshen2019@foxmail.com}\\
\small School of Mathematics and Statistics, and\\
\small Key Laboratory of Nonlinear Analysis and Applications (Ministry of Education),\\ \small Central China Normal University, Wuhan 430079, P. R. China}



\date{}
\maketitle
\noindent {\bf Abstract}:
Let  $P_{10}$ be a path on $10$ vertices. A graph is said to be $P_{10}$-free if it does not contain $P_{10}$ as an induced subgraph. The well-known Erd\H{o}s-Gy\'{a}rf\'{a}s Conjecture states that every graph with minimum degree at least three has a cycle whose length is a power of $2$. In this paper, we show that every $P_{10}$-free graph with minimum degree at least three contains a cycle of length $4$ or $8$. This implies that the conjecture is true for $P_{10}$-free graphs.

\vskip 0.1cm


\section{Introduction}
All graphs considered here are finite and simple. Let $G$ be a graph. The vertex set, the edge set and the minimum degree of $G$ are denoted by $V(G)$, $E(G)$ and $\delta(G)$, respectively.
For a vertex $v\in V(G)$, we denote by $N_{G}(v)$ the neighbors of $v$ in $G$. For $S \subseteq V(G)$, let $N_G(S)=\cup_{x \in S} N_G(x) - S$.
For convenience, we write $N(v)$ and $N(S)$ for $N_G(v)$ and $N_G(S)$, respectively.
Denote by $G[S]$  the  subgraph of $G$ induced by $S$.
For $X, Y\subseteq V(G)$, $E_G(X,Y)$ represents the set consisting
of all edges in $G$ with one end in $X$ and the other in $Y$.
Let $H$ be a graph. We say that $G$ is \emph{$H$-free} if it does not contain  an induced subgraph isomorphic to $H$.
For $H\subseteq G$, we use $G-H$ to denote the subgraph of $G$ induced by $V(G)-V(H)$.
Let $x,y$ be two distinct vertices of $G$. An \emph{$(x,y)$-path} in $G$ is a path from $x$ to $y$. The length of a shortest $(x,y)$-path in $G$ is denoted by $d_G(x,y)$. If $P$ is an $(x,y)$-path in $G-E(H)$ such that $V(P)\cap V(H)=\{x,y\}$, then we call $P$ an \emph{$H$-path} in $G$. (To specify the end vertices of $P$, we also call $P$ an \emph{$H$-$(x,y)$-path}).
We denote a path on $k$ vertices (resp. a cycle on $k$ vertices) by $P_k$ (resp. $C_k$). A cycle of length three is also called a triangle. An edge is a \emph{triangulated edge} if it lies in a triangle. The length of a cycle $C$ is denoted by $\ell(C)$. A \emph{chord} of a cycle $C$ is an edge of $G-E(C)$ joining  two vertices of $C$. A \emph{hole} of $G$ is an induced cycle, that is, a cycle without a chord. An $m$-hole is a hole of length $m$.

The well-known Erd\H{o}s-Gy\'{a}rf\'{a}s Conjecture \cite{10} states that every graph with minimum degree at least three has a cycle whose length is a power of $2$. Markstr\"{o}m \cite{02} has shown that any cubic counterexample to this conjecture must have at least $30$ vertices. Nowbandegani and Esfandiari \cite{03} showed that any bipartite counterexample must have at least $32$ vertices.  Nowbandegani et al. \cite{01} proved that any cubic claw-free counterexample must have at least $114$ vertices. Moreover, the conjecture was confirmed in the following  graph classes: $K_{1,m}$-free graphs with minimum degree at least $m+1$ or maximum degree at least $2m-1$ \cite{04}, planar claw-free graphs \cite{05}, $3$-connected cubic planar graphs \cite{06}, $P_8$-free graphs \cite{09} and some Cayley graphs \cite{07,08}. In \cite{01}, Nowbandegani et al. showed that every claw-free graph with minimum degree at least  $3$ has a cycle whose length is $2^k$ or $3 \cdot 2^k$, for some positive integer $k$. Very recently, Liu and Montgomery \cite{11} proved, using a new technique for constructing even cycles, that there exists a constant $c$ such that every graph with average degree at least $c$ contains a cycle whose length is a power of $2$.


In this paper, we confirm the Erd\H{o}s-Gy\'{a}rf\'{a}s Conjecture for $P_{10}$-free graphs by providing the following theorem.

\begin{thm}\label{thm1.1}
Every $P_{10}$-free graph with minimum degree at least three contains a $C_4$ or $C_8$.
\end{thm}

\section{Preliminaries}
In this section, we  give some notations and lemmas used in this article.
For two integers $s$ and $t$ with $s\leq t$ , we use $[s,t]$ to denote the set of all integers between $s$ and $t$.

For a cycle $C=x_1x_2\cdots x_mx_1$,
we specify that the positive direction of $C$ is the direction in which the subscripts of vertices increase successively, and the opposite is the negative direction. For $i,j\in[1,m]$, $x_i\overrightarrow{C}x_j$ is a path from $x_i$ to $x_j$ along $C$ with a positive direction. $x_i\overleftarrow{C}x_j$ is a path from $x_i$ to $x_j$ along $C$ with a negative direction. For convenience, we let
$$
N_C(x):=N(x)\cap V(C),~\forall~ x\in V(G).
$$
Set
    $$
    A_i:=N_{G-C}(x_i), ~1\leq i\leq m.
    $$
Moreover, we let
$$
I_C:=\{x_i:~i\in [1,m]~\makebox{and}~A_i\cap A_{i+1}\neq\emptyset\}
$$
and
$$
I_C^+:=\{x_i:~x_{i-1}\in I_C\},
$$
where the subscripts of $x_i$ and $A_i$ are taken modulo $m$. Clearly, if $C$ is an $m$-hole with $m\geq 4$, then $x_i\in I_C$ if and only if the edge $x_ix_{i+1}$ lies in a triangle.

\begin{defi}\label{def-2-1}
A hole $C$ of length at least $5$ in $G$ is good if
\begin{itemize}
  \item there exists no hole $C'$ in $G$ such that $5\leq \ell(C')<\ell(C)$, and
  \item  subject to this, $|I_C|$ is as large as possible.
\end{itemize}
If, in addition, $\ell(C)=m$, then we call $C$  a good $m$-hole.
\end{defi}
\begin{defi}\label{def-2-2}
Let $H$ be a subgraph of a graph $G$ and let $X\subseteq V(H)$. If $P:=u_1\ldots u_t$ is an induced path of $G-H$ such that
$$
E_G(V(P),V(H))=\{u_1x:~x\in X\},
$$
then we call $P$ a good path for $(H, X)$. If, in addition, $X=\{x\}$, we call $xu_1\ldots u_t$ a good $(H,x)$-path.
\end{defi}

Now, we  show some technical lemmas for graphs without $C_4$, which will be used in the proof of Theorem \ref{thm1.1}.

\begin{lem}\label{lem-2-1}
Let $G$ be a graph with $\delta(G) \geq 3$ and $C$ a cycle of length at least $4$ in $G$. If $G$ does not contain $C_4$, then $G[V(C)]$ has an $m$-hole for some integer $m$ with $5\leq m\leq \ell(C)$.
\end{lem}

\noindent{\bf Proof}\ Pick a cycle $D$ in $G[V(C)]$ such that

(i) $\ell(D)\geq 4$, and

(ii) subject to (i), $\ell(D)$ is as small as possible.

\noindent Then, $\ell(D)\leq \ell(C)$. As $G$ does not contain $C_4$, $\ell(D)\geq 5$.
If $D$ has no chord, then $D$ is the desired $m$-hole with $m=\ell(D)$. By way of contradiction, assume that $D$ contains a chord $xy$. Set $C'=x\overrightarrow{D}yx$ and
$C''=x\overleftarrow{D}yx$. Then, $C'$ and $C''$ are cycles  in $G$ such that $\ell(C')+\ell(C'')=\ell(D)+2$. By symmetry, we may assume that $\ell(C')\geq \ell(C'')$.
Then, $C'$ is a cycle in $G[V(C)]$ with $4\leq \ell(C')\leq\ell(D)-1$, contrary to the choice of $D$. Hence, Lemma \ref{lem-2-1} is true. \qed

Note that every graph with minimum degree at least $3$ admits a cycle of length at least $4$. By Lemma \ref{lem-2-1}, the following lemma, due to Nowbandegani et al., holds.

\begin{lem}\label{lem-2-2}\emph{\cite{01}}
Let $G$ be a graph with $\delta(G) \geq 3$. If $G$ does not contain $C_4$ as a subgraph, then $G$ has an $m$-hole for some $m \geq 5$.
\end{lem}

\begin{lem}\label{lem-2-3}
Let $G$ be a graph and $u,v,v'$ three vertices of $G$ such that $v,v'\in N(u)$. Let $A$ be a subset of $V(G)-\{u,v,v'\}$ such that
$$
\min~\{|N(v)\cap A|,|N(v')\cap A|\}>0.
 $$
 If $G$ does not contain $C_4$ as a subgraph, then there exist two independent edges from $\{v,v'\}$ to $A$.
\end{lem}
\noindent{\bf Proof}\ By way of contradiction, assume that Lemma \ref{lem-2-3} is false.
Then, $v$ and $v'$ share a common neighbor, say $w$, in $A$.  It follows that $vwv'uv$ is a $C_4$ in $G$, a contradiction.
Hence, Lemma \ref{lem-2-3} is true. \qed

\begin{lem}\label{lem-2-4}
Let $G$ be a graph and let $C:=x_1x_2\ldots x_mx_1$ be a cycle
 in $G$ with $5\leq m\leq 7$. If $G$  contains neither $C_4$ nor $C_8$ as a subgraph, then $|I_C|\leq 7-m$.
\end{lem}

\noindent{\bf Proof}\ By way of contradiction, assume that Lemma \ref{lem-2-4} is false. Then, $t:=|I_C|\geq 8-m$. Denote $I_C:=\{x_{i_1},\ldots,x_{i_t}\}$, where $1\leq i_1<\ldots<i_t\leq m$. For each $j\in [1,t]$, let $y_{i_j}\in A_{i_j}\cap A_{i_j+1}$, that is, $y_{i_j}\in N_{G-C}(x_{i_j})\cap N_{G-C}(x_{i_j+1})$. We claim that
\begin{equation}\label{eqn-2-1}
    N(y_i)\cap \{x_i,x_{i+1},x_{i+2},x_{i+3}\}=\{x_i,x_{i+1}\},~\forall ~x_i\in I_C.
\end{equation}
By  way of contradiction, assume that (\ref{eqn-2-1}) is false. Then,  $x_j\in N(y_i)$ holds for some $j\in\{i+2,i+3\}$. Noting that $x_{j-2}\in N(y_i)$, we conclude that $x_{j-2}x_{j-1}x_jy_ix_{j-2}$ is a  $C_4$ in $G$, a contradiction. Hence, (\ref{eqn-2-1}) is true. By symmetry, we also have
$$
    N(y_i)\cap \{x_{i-2},x_{i-1},x_i,x_{i+1}\}=\{x_i,x_{i+1}\},~\forall ~x_i\in I_C.
$$
This together with (\ref{eqn-2-1}) and $m\leq 7$ implies that
$$
 \{x_i,x_{i+1}\}\subseteq   N_C(y_i)\subseteq\{x_i,x_{i+1}\}\cup(\{x_{i+4}\}\cap\{x_{i-3}\}),~\forall ~x_i\in I_C.
$$
Hence, $y_{i_1},\ldots,y_{i_t}$ are distinct vertices of $V(G)-V(C)$. Let $C'$ be the cycle obtained from $C$ by replacing $\cup_{j=1}^{8-m}\{x_{i_j}x_{i_j+1}\}$ with
$\cup_{j=1}^{8-m}\{x_{i_j}y_{i_j}x_{i_j+1}\}$. Then, $C'$ is a  $C_8$ in $G$, a contradiction. Hence, Lemma \ref{lem-2-4} is true.\qed

\begin{lem}\label{lem-2-5}
Let $G$ be a graph with $\delta(G) \geq 3$ and let $C:=x_1x_2\ldots x_mx_1$ be a good hole in $G$. If $G$  contains neither $C_4$ nor $C_8$ as a subgraph, then for each $i\in [1,m]$, there exists  a good path  for $(C,X_i)$ with order  $\min~\{\lfloor m/ 2\rfloor-1, ~2\}$, where
$$
X_i:=\left\{\begin{array}{ll}
       \{x_i,x_{i+1}\} & \mbox{if $x_i\in I_C$} \cr
       \{x_{i-1},x_i\} & \mbox{if $x_i\in I_C^+\setminus I_C$}\cr
       \{x_i\} & \mbox{if $x_i\notin I_C\cup I_C^+$.}
     \end{array}
     \right.
$$
\end{lem}

\noindent{\bf Proof}\ As $C$ is a hole, $|A_i|=|N(x_i)-\{x_{i-1},x_{i+1}\}|\geq 1$.  We claim that
\begin{equation}\label{eqn-2-2}
N_C(u)\subseteq\{x_{i-1}, x_i, x_{i+1}\},~\forall ~u\in A_i.
\end{equation}
By way of contradiction, assume that (\ref{eqn-2-2}) is false. Then, there exists $u\in A_i$ such that $N(u)\cap (V(C)-\{x_{i-1}, x_i, x_{i+1}\})$ contains at least one vertex, say $x_j$.  Set $C'=x_i\overrightarrow{C}x_jux_i$ and
$C''=x_i\overleftarrow{C}x_jux_i$. Then, $C'$ and $C''$ are cycles of length at least $4$ in $G$ such that $\ell(C')+\ell(C'')=\ell(C)+4$. By symmetry, we may assume that $\ell(C')\leq \ell(C'')$. Then,
$$
4\leq\ell(C')\leq \frac{\ell(C)+4}{2}<\ell(C).
$$
By Lemma \ref{lem-2-1}, $G$ has a hole $C^*$ such that
$5\leq \ell(C^*)\leq\ell(C')<\ell(C)$, a contradiction. Hence, (\ref{eqn-2-2}) is true.

Choose $u_1\in A_i$ such that
\begin{itemize}
  \item $u_1\in A_i\cap A_{i+1}$ if $x_i\in I_C$, and
  \item $u_1\in A_{i-1}\cap A_i$ if $x_i\in I^+_C\setminus I_C$.
\end{itemize}
If $\{x_{i-1},x_{i+1}\}\subseteq N(u_1)$, then $x_{i-1}x_ix_{i+1}u_1x_{i-1}$ is a $C_4$ in $G$, a contradiction. Hence, $\{x_{i-1},x_{i+1}\}\not\subseteq N(u_1)$. This together with (\ref{eqn-2-2}) and the choice of $u_1$ implies that
\begin{equation}\label{eqn-2-3}
N_C(u_1)=X_i.
\end{equation}
It follows from (\ref{eqn-2-3}) that $u_1$ is a good path for $(C,X_i)$, and hence Lemma \ref{lem-2-5} holds for $m=5$. In the following, we assume that $m\geq 6$.

By (\ref{eqn-2-3}), $|N(u_1)-V(C)|\geq d(u_1)-2\geq 1$, and hence $N(u_1)-V(C)\neq\emptyset$. Choose $u_2\in N(u_1)-V(C)$ such that

(i) $|N(u_2)\cap \{x_i\}|$ achieves the minimum, and

(ii) subject to (i), $|N_C(u_2)|$ is as small as possible.

\noindent We claim that
\begin{equation}\label{eqn-2-4}
x_i\notin N(u_2).
\end{equation}
By way of contradiction, assume that $x_i\in N(u_2)$.  Then, by the choice of $u_2$, we have  $N(u_1)-V(C)\subseteq N(x_i)$. If $N(u_1)-V(C)$ contains a vertex $u_2'\neq u_2$, then $u_2'\in N(x_i)$, and hence $u_2u_1u_2'x_iu_2$ is a $C_4$ in $G$, a contradiction. Hence, $N(u_1)-V(C)=\{u_2\}$. As $\delta(G)\geq 3$, $|N_C(u_1)|\geq 2$.
This together with (\ref{eqn-2-3})   implies that $N_C(u_1)=\{x_i, x_j\}$ holds for some $x_j$ with $j\equiv i\pm 1~(\makebox{mod}~ m)$. Noting that $x_ix_j\in E(C)$, we conclude that $x_iu_2u_1x_jx_i$ is a $C_4$ in $G$, a contradiction. Hence, (\ref{eqn-2-4}) is true.

If $N_C(u_2)=\emptyset$, then by (\ref{eqn-2-3}), we see that $u_1u_2$ is a good path for $(C,X_i)$, which implies that Lemma \ref{lem-2-5} is true. By way of contradiction, assume that Lemma \ref{lem-2-5} is false, then $N_C(u_2)\neq\emptyset$. Say $x_j\in N_C(u_2)$ for some $j\in [1, m]$. By (\ref{eqn-2-4}), $j\neq i$. Set $D'=x_i\overrightarrow{C}x_ju_2u_1x_i$ and
$D''=x_i\overleftarrow{C}x_ju_2u_1x_i$. Then, $D'$ and $D''$ are cycles of length at least $4$ in $G$ such that $\ell(D')+\ell(D'')=\ell(C)+6$. By symmetry, we may assume that $\ell(D')\leq \ell(D'')$. Then,
$$
4\leq\ell(D')\leq \frac{\ell(C)+6}{2}.
$$
By Lemma \ref{lem-2-1}, $G$ has a hole $C^*$ with
$5\leq \ell(C^*)\leq\ell(D')\leq \frac{\ell(C)+6}{2}\leq \ell(C)$. As $C$ is a good $m$-hole, $\ell(C^*)\geq\ell(C)$, and hence $\ell(C^*)=\ell(D')=\ell(C)=m=6$, which in turn means  $j=i+3~(\makebox{mod $m$})$. It follows that both $D'$ and $D''$ are $6$-holes. Thus,
\begin{equation}\label{eqn-2-5}
\makebox{$N_C(u_1)=\{x_i\}$ and $N_C(u_2)=\{x_{i+3}\}$}.
\end{equation}
As $\delta(G)\geq 3$, $u_1$ has a neighbor $u_2'\in V(G)-(V(C)\cup\{u_2\})$. If $x_i\in N(u_2')$, then $x_iu_1u_2'x_i$ is a triangle in $G$, and hence $u_1\in I_{D'}\cap I_{D''}$ (see Figure 1).
\begin{figure}\label{fig-1}
\psfrag{A}{$x_i$}
\psfrag{B}{$x_{i+1}$}
\psfrag{C}{$x_{i+2}$}
\psfrag{D}{$x_{i+3}$}
\psfrag{E}{$x_{i+4}$}
\psfrag{F}{$x_{i+5}$}
\psfrag{I}{$u_1$}
\psfrag{J}{$u_2$}
\psfrag{K}{$u_2'$}
      \begin{center}
     \includegraphics[width=55mm]{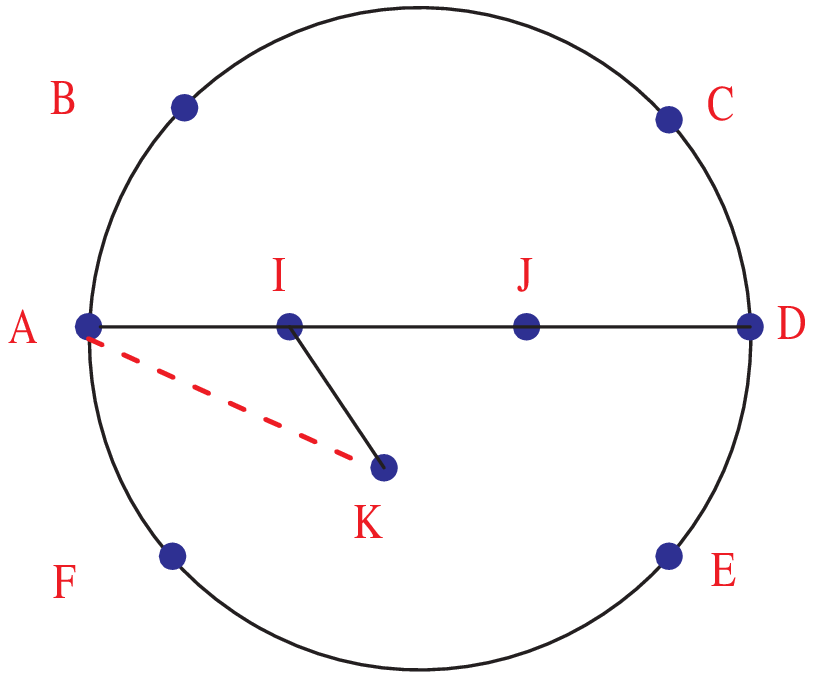}\\
   \caption{Illustration of the nonexistence of the edge $x_iu_2'$.}
   \end{center}
\end{figure}
This together with the choice of $C$ implies that $I_C\neq\emptyset$. Say $x_k\in I_C$ for some $k\in [1,6]$.
\begin{itemize}
  \item If $x_k\in \{x_i,x_{i+1},x_{i+2}\}$, then $D'$ is a $6$-hole with $I_{D'}\supseteq\{u_1,x_k\}$;
  \item if $x_k\in\{x_{i+3},x_{i+4},x_{i+5}\}$, then $D''$ is a $6$-hole with $I_{D''}\supseteq\{u_1,x_k\}$.
\end{itemize}
In each case, we get  a contradiction to Lemma \ref{lem-2-4}. Hence, $x_i\notin N(u_2')$. This together with the choice of $u_2$ implies that $N_C(u_2')\neq\emptyset$. By an argument similar to that in the proof of (\ref{eqn-2-5}), we can derive that $N_C(u_2')=\{x_{i+3}\}$. This together with (\ref{eqn-2-5}) implies that $u_1u_2x_{i+3}u_2'u_1$ is a $C_4$ in $G$, a contradiction. Hence, Lemma \ref{lem-2-5}
is true. \qed

\vspace{3mm}

For $m\geq 5$,  let $t_m(G)$ be  the minimum non-negative integer $t$ such that every $m$-hole of $G$ admits at most $t$  triangulated edges. (If $G$ does not contain any $m$-hole, we let $t_m(G)=0$).

A theta-graph \emph{$\theta(a,b,c)$}, $1\leq a\leq b\leq c$, $b\geq 2$,  is a simple graph consisting of $3$ internally disjoint paths of lengths $a$, $b$ and $c$ between a pair of vertices of degree $3$. We conclude this section with the following lemma on theta-graphs.

\begin{lem}\label{lem-2-6}
Let $G$ be a $P_{10}$-free graph with minimum degree at least $3$ and let $H$ be a subgraph of $G$ isomorphic to $\theta(2,3,3)$. If $t_5(G)=0$, then $G$ admits a $C_4$ or $C_8$.
\end{lem}

\begin{figure}\label{fig-2}
\psfrag{Y}{$y$}
\psfrag{A}{$x_1$}
\psfrag{B}{$x_2$}
\psfrag{C}{$x_3$}
\psfrag{D}{$x_4$}
\psfrag{E}{$x_5$}
\psfrag{F}{$x_6$}
      \begin{center}
     \includegraphics[width=55mm]{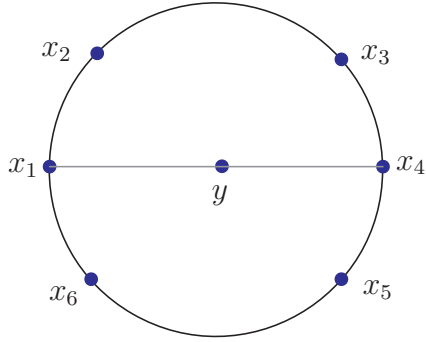}\\
   \caption{The graph $\theta(2,3,3)$.}
   \end{center}
\end{figure}

\noindent{\bf Proof}\ By way of contradiction, assume that $G$ contains neither $C_4$ nor $C_8$.
Mark the vertices of $H$ as in Figure 2.  Set $C':=x_1x_2x_3x_4yx_1$ and
$C'':=x_4x_5x_6x_1yx_4$.

As $G$ does not contain $C_4$, both $C'$ and $C''$ are $5$-holes of $G$. If $H$ is not any induced subgraph of $G$, then $x_ix_j\in E(G)$ holds for some integers $i$, $j$ with  $i\in\{2,3\}$ and $j\in\{5,6\}$. By symmetry, we may assume $i=2$. If $j=5$, then $x_2x_3x_4x_5x_2$ is a $C_4$ in $G$; if $j=6$, then  $x_1x_2$ is a triangulated edge in the $5$-hole $C'$, which implies $t_5(G)\geq 1$. Either way leads a contradiction. Hence, $H$ is an induced subgraph of $G$.

Set
$$
 Y_i:=N(x_i)-V(H), ~i=2,3,5,6.
$$
As $\delta(G)\geq 3$, $Y_i\neq\emptyset$, $i=2,3,5,6$.
Let  $y_2\in Y_2$. If $x_1\in N(y_2)$, then $x_1x_2$ is a triangulated edge in the $5$-hole $C'$, a contradiction. Hence, $x_1\notin N(y_2)$.  Likewise, $x_3\notin N(y_2)$.  If $x_5\in N(y_2)$, then $x_2y_2x_5x_6x_1yx_4x_3x_2$ is a $C_8$ in $G$, a contradiction. Hence, $x_5\notin N(y_2)$. Finally,  $x_4, x_6,y\notin N(y_2)$, since otherwise $x_4y_2x_2x_3x_4$, $x_6y_2x_2x_1x_6$ or $yy_2x_2x_1y$ is a $C_4$ in $G$, a contradiction. Therefore,
    $N(y_2)\cap V(H)=\{x_2\}$. Likewise, we have
    \begin{equation}\label{eqn-2-6}
N(Y_i)\cap V(H)=\{x_i\}, ~\forall~i\in \{2,3,5,6\}.
\end{equation}
If there exists two $H$-paths $P$ and $Q$ of length $3$ connecting $x_2, x_5$ and $x_3,x_6$ respectively, say $P:=x_2y_2y_5x_5$ and $Q:=x_3y_3y_6x_6$, then by (\ref{eqn-2-6}), we can derive  that
$$
N(y_i)\cap V(H)=\{x_i\},~i=2,3,5,6.
$$
It follows that $N(y_i)\cap V(H)$, $i=2,3,5,6$, are distinct subsets of $V(H)$, and hence
$y_2,y_3,y_5,y_6$ are distinct vertices of $G-H$. Thus,  $x_2y_2y_5x_5x_6y_6y_3x_3x_2$ is a $C_8$ in $G$ (see Figure 3), a contradiction.
By symmetry, we may assume that there exists no $H$-$(x_2,x_5)$-path of length $3$ in $G$.
\begin{figure}\label{fig-3}
\psfrag{y}{$y$}
\psfrag{P}{$P$}
\psfrag{Q}{$Q$}
\psfrag{A}{$x_1$}
\psfrag{B}{$x_2$}
\psfrag{C}{$x_3$}
\psfrag{D}{$x_4$}
\psfrag{E}{$x_5$}
\psfrag{F}{$x_6$}
\psfrag{G}{$y_2$}
\psfrag{H}{$y_5$}
\psfrag{I}{$y_6$}
\psfrag{J}{$y_3$}
\begin{center}
     \includegraphics[width=55mm]{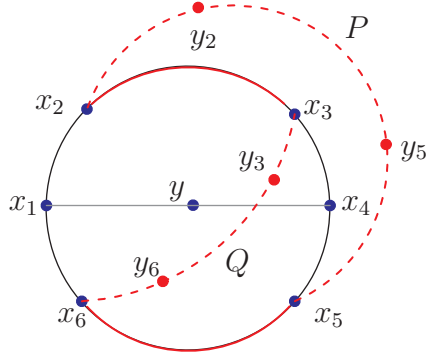}\\
   \caption{Illustration of the nonexistence of the two $H$-paths $P$ and $Q$.}
\end{center}
\end{figure}

If there exists an  $H$-$(x_2,x_5)$-path $R$ of length at most $5$ in $G$, then $\ell(R)\in\{1,2,4,5\}$, and hence
$$
C^*:=\left\{\begin{array}{ll}
x_2x_5x_4x_3x_2,&\makebox{if $\ell(R)=1$}\cr
x_2\overrightarrow{R}x_5x_6x_1yx_4x_3x_2,&\makebox{if $\ell(R)=2$}\cr
x_2\overrightarrow{R}x_5x_4yx_1x_2,&\makebox{if $\ell(R)=4$}\cr
x_2\overrightarrow{R}x_5x_6x_1x_2,&\makebox{if $\ell(R)=5$}
\end{array}\right.
$$
is a $C_4$ or $C_8$ in $G$, a contradiction. Thus,
\begin{equation}\label{eqn-2-7}
\makebox{there exists no $H$-$(x_2,x_5)$-path of length at most $5$ in $G$.}
\end{equation}

Now, we prove the following claim.
\begin{claim} \label{cla-2-1}
For $i=2,5$, $G$ contains a good  $(H,x_i)$-path of length $3$.
\end{claim}
\noindent{\bf Proof}\
By symmetry, it suffices to show Claim \ref{cla-2-1} for $i=2$.
Consider any vertex $y_2$ in $Y_2$. By (\ref{eqn-2-6}), $N(y_2)\cap V(H)=\{x_2\}$.
As $\delta(G)\geq 3$, $y_2$ has at least two neighbors, say $z_2$ and $z_2'$, in $V(G)-V(H)$.   We claim that
\begin{equation}\label{eqn-2-8}
\min~\{|N(z_2)\cap\{x_2,x_4\}|,|N(z_2')\cap\{x_2,x_4\}|\}=0.
\end{equation}
By way of contradiction, assume  (\ref{eqn-2-8}) is false. Then, both $z_2$ and $z_2'$ has a neighbor in $\{x_2,x_4\}$. By applying Lemma \ref{lem-2-3} with $(u,v,v',A):=(y_2,z_2,z_2',\{x_2,x_4\})$, we can derive that $E_G(\{z_2,z_2'\},\{x_2,x_4\})$ contains two independent edges.
By renaming $z_2$ and $z_2'$ (if necessary), we may assume that $x_2\in N(z_2)$ and $x_4\in N(z_2')$. Then, $x_2z_2y_2z_2'x_4x_5x_6x_1x_2$ is a $C_8$ in $G$, a contradiction. Hence, (\ref{eqn-2-8}) is true.

By symmetry, we may assume that $N(z_2)\cap\{x_2,x_4\}=\emptyset$.
We claim that
\begin{equation}\label{eqn-2-9}
N(z_2)\cap V(H)=\emptyset.
\end{equation}
Otherwise,  we have $N(z_2)\cap \{y,x_1,x_3,x_5,x_6\}\neq\emptyset$.
\begin{itemize}
   \item If $y\in N(z_2)$, then $x_2y_2z_2yx_4x_5x_6x_1x_2$ is a $C_8$ in $G$;
   \item if $x_k\in N(z_2)$ for some $k\in\{1,3\}$, then $x_2y_2z_2x_kx_2$ is a $C_4$ in $G$;
   \item if $x_5\in N(z_2)$, then $x_2y_2z_2x_5$ is an $H$-path in $G$  contradicting (\ref{eqn-2-7});
   \item if $x_6\in N(z_2)$, then $x_2y_2z_2x_6x_5x_4yx_1x_2$ is a $C_8$ in $G$.
\end{itemize}
Thus in all cases, we obtain a contradiction. Hence, (\ref{eqn-2-9}) is true.

\begin{figure}\label{fig-4}
\psfrag{y}{$y$}
\psfrag{A}{$x_1$}
\psfrag{B}{$x_2$}
\psfrag{C}{$x_3$}
\psfrag{D}{$x_4$}
\psfrag{E}{$x_5$}
\psfrag{F}{$x_6$}
\psfrag{G}{$y_2$}
\psfrag{H}{$z_2$}
\psfrag{I}{$z_2'$}
\psfrag{J}{$u_2$}
\psfrag{K}{$u_2'$}
      \begin{center}
     \includegraphics[width=60mm]{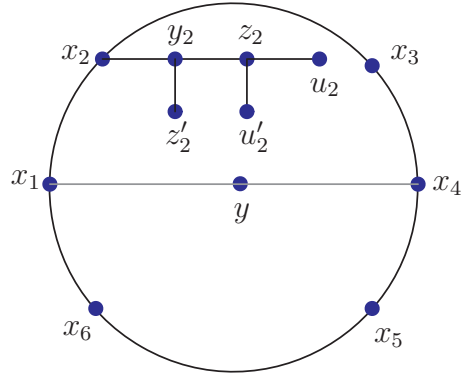}\\
   \caption{Vertices used in the proof of Claim \ref{cla-2-1}.}
   \end{center}
\end{figure}

It follows from (\ref{eqn-2-9}) that $N(z_2)\cap (V(H)\cup\{y_2\})=\{y_2\}$. As $\delta(G)\geq 3$, $z_2$ has two neighbors, say $u_2$ and $u_2'$,  in $V(G)-(V(H)\cup\{y_2\})$ (see Figure 4).  We claim that
\begin{equation}\label{eqn-2-10}
\min~\{|N(u_2)\cap\{y_2,y\}|,|N(u_2')\cap\{y_2,y\}|\}=0.
\end{equation}
By way of contradiction, assume  (\ref{eqn-2-10}) is false.  Then, both $u_2$ and $u_2'$ has a neighbor in $\{y_2,y\}$. By applying lemma \ref{lem-2-3} with $(u,v,v',A):=(z_2,u_2,u_2',\{y_2,y\})$, we can derive that $E_G(\{u_2,u_2'\},\{y_2,y\})$ contains two independent edges. By renaming $u_2$ and $u_2'$ (if necessary), we may assume that $y_2\in N(u_2)$ and $y\in N(u_2')$. Then, $x_2y_2u_2z_2u_2'yx_4x_3x_2$ is a $C_8$ in $G$, a contradiction. Hence, (\ref{eqn-2-10}) is true.

By symmetry, we may assume that
\begin{equation}\label{eqn-2-11}
N(u_2)\cap\{y_2,y\}=\emptyset.
\end{equation}
If $N(u_2)\cap V(H)\neq\emptyset$, then by (\ref{eqn-2-11}),  we can derive that $u_2x_k\in E(G)$  for some $k\in [1,6]$. Set
$$
D:=\left\{
  \begin{array}{ll}
   x_2y_2z_2u_2x_1yx_4x_3x_2, & \hbox{if $k=1$}\cr
   x_2y_2z_2u_2x_2, & \hbox{if $k=2$}\cr
   x_2y_2z_2u_2x_3x_4yx_1x_2, & \hbox{if $k=3$}\cr
x_2y_2z_2u_2x_4x_5x_6x_1x_2, & \hbox{if $k=4$}\cr
x_2y_2z_2u_2x_5x_4yx_1x_2, & \hbox{if $k=5$}\cr
x_2y_2z_2u_2x_6x_5x_4x_3x_2,& \hbox{if $k=6$.}
\end{array}
\right.
$$
Then, $D$ is a $C_4$ or $C_8$ in $G$, a contradiction. Hence,
\begin{equation}\label{eqn-2-12}
N(u_2)\cap V(H)=\emptyset.
\end{equation}

Recall that $N(y_2)\cap V(H)=\{x_2\}$. By (\ref{eqn-2-9}), (\ref{eqn-2-11}) and (\ref{eqn-2-12}), we see that $x_2y_2z_2u_2$ is a good $(H,x_2)$-path in $G$.
Likewise, $G$ contains a good $(H,x_5)$-path of length $3$. Hence, Claim \ref{cla-2-1} is true. \qed

It follows from  Claim \ref{cla-2-1} that for each $i=2,5$,  $G$ contains a good good  $(H,x_i)$-path, say $x_iy_iz_iu_i$,  of length $3$. Set $R_2:=y_2z_2u_2$ and $R_5:=y_5z_5$. Then, both $R_2$ and $R_5$ are induced paths of $G-V(H)$, and
\begin{equation}\label{eqn-2-13}
E_G(V(R_i),V(H))=\{y_ix_i\}, ~i=2,5.
\end{equation}

If $V(R_2)\cap V(R_5)\neq\emptyset$, then $G[V(R_2)\cup V(R_5)]$ contains a $(y_2,y_5)$-path $Q$ of length at most $3$. It follows that $x_2y_2Qy_5x_5$ is an $H$-$(x_2,x_5)$-path of length at most $5$ in $G$, contrary to (\ref{eqn-2-7}).
     Hence, $V(R_2)\cap V(R_5)=\emptyset$.  Similarly, we have
 \begin{equation}\label{eqn-2-14}
E_G(V(R_2),V(R_5))\subseteq\{u_2z_5\}.
\end{equation}
Recall that $y_2z_2u_2$ and $y_5z_5u_5$ are induced paths of $G-V(H)$. If $u_2z_5\notin E(G)$, then by (\ref{eqn-2-13})  and (\ref{eqn-2-14}), we can derive  that $u_2z_2y_2x_2x_1yx_4x_5y_5z_5$ is an induced $P_{10}$ in $G$, a contradiction. Hence,
\begin{equation}\label{eqn-2-15}
u_2z_5\in E(G).
\end{equation}

For $i=3,6$, let $y_i\in Y_i$. It follows from
(\ref{eqn-2-6}) that
 \begin{equation}\label{eqn-2-16}
N(y_i)\cap V(H)=\{x_i\},~i=3,6.
\end{equation}
This together with (\ref{eqn-2-13}) implies that $y_3\neq y_6$ and
$\{y_3,y_6\}\cap\{y_2,z_2,u_2,y_5,z_5\}=\emptyset$. If $\{y_3,y_6\}\subseteq N(u_2)$, then $x_3y_3u_2y_6x_6x_1yx_4x_3$ is a $C_8$ in $G$, a contradiction.
\begin{figure}\label{fig-5}
\psfrag{y}{$y$}
\psfrag{A}{$x_1$}
\psfrag{B}{$x_2$}
\psfrag{C}{$x_3$}
\psfrag{D}{$x_4$}
\psfrag{E}{$x_5$}
\psfrag{F}{$x_6$}
\psfrag{G}{$y_2$}
\psfrag{H}{$z_2$}
\psfrag{I}{$u_2$}
\psfrag{J}{$y_5$}
\psfrag{K}{$z_5$}
\psfrag{L}{$y_3$}
\psfrag{M}{$y_6$}
      \begin{center}
     \includegraphics[width=60mm]{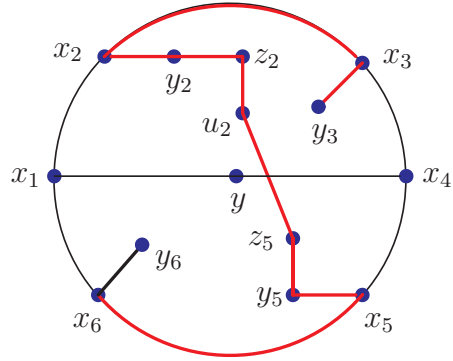}\\
   \caption{An induced $P_{10}$ in case $N(y_3)\cap\{y_2,z_2,y_5,z_5\}=\emptyset$.}
   \end{center}
\end{figure}
Hence, $\{y_3,y_6\}\not\subseteq N(u_2)$. By symmetry, we may assume that $y_3\notin N(u_2)$. If $N(y_3)\cap\{y_2,z_2,y_5,z_5\}=\emptyset$, then by (\ref{eqn-2-13}), (\ref{eqn-2-14}), (\ref{eqn-2-15})  and (\ref{eqn-2-16}), we can derive that
$y_3x_3x_2y_2z_2u_2z_5y_5x_5x_6$ is an induced $P_{10}$ in $G$ (see Figure 5),
a contradiction. Hence, $N(y_3)\cap\{y_2,z_2,y_5,z_5\}\neq\emptyset$. Set
$$
D^*=\left\{
     \begin{array}{ll}
       x_3y_3y_2x_2x_3, & \hbox{if $y_2\in N(y_3)$} \cr
       x_3y_3z_2y_2x_2x_1yx_4x_3, & \hbox{if $z_2\in N(y_3)$} \cr
       x_3y_3z_5y_5x_5x_6x_1x_2x_3, & \hbox{if $z_5\in N(y_3)$} \cr
x_3y_3y_5x_5x_4yx_1x_2x_3, & \hbox{if $y_5\in N(y_3)$}.
     \end{array}
   \right.
$$
Then, $D^*$ is a  $C_4$ or $C_8$ in $G$, a contradiction. This completes the proof of Lemma \ref{lem-2-6}. \qed

\section{\bf Proof of Theorem \ref{thm1.1}.}
\noindent  Suppose, by contradiction, that $G$ contains neither $C_4$ nor $C_8$. By Lemma \ref{lem-2-2}, $G$ contains a hole of length at least $5$. Let $C:=x_1x_2\cdots x_mx_1$ be a hole in $G$ such that
\begin{wst}
    \item[(T1)] $\ell(C)\geq 5$;
    \item[(T2)] subject to (T1), $\ell(C)$ is as small as possible;
    \item[(T3)] subject to (T1) and (T2), $|I_C|$ is as large as possible.
\end{wst}
Then, $C$ is a good $m$-hole. Define $A_i$, $i\in [1,m]$, as that in Section 2. We show some claims, the first one of which is used frequently.

\begin{claim} \label{claim-3-1}
Let $D$ be a cycle of $G$. If $4\leq\ell(D)\leq\ell(C)$, then $D$ is an $m$-hole and $|I_D|\leq |I_C|$. As a consequence, there exists no cycle of length $k$ in $G$ with $4\leq k<m$.
\end{claim}
\noindent{\bf Proof}\ By Lemma \ref{lem-2-1}, $G[V(D)]$ contains a hole $D^*$ with $5\leq \ell(D^*)\leq\ell(D)\leq\ell(C)$. This together with the choice of $C$ implies that $D^*$ is an $m$-hole, and hence $D=D^*$. By the choice of $C$, we see that $|I_D|\leq |I_C|$. Hence, Claim \ref{claim-3-1} is true. \qed

\begin{claim} \label{claim-3-2}
$\ell(C)\leq 6$.
\end{claim}

\noindent{\bf Proof}\ By way of contradiction, assume that Claim \ref{claim-3-2} is false. Then, $m\geq 7$. By Lemma \ref{lem-2-5}, for each $i\in [1,m]$, $G$ has a good path, say  $y_iz_i$,  for $(C,X_i)$ with order $2$, where
$$
X_i:=\left\{\begin{array}{ll}
       \{x_i,x_{i+1}\} & \mbox{if $x_i\in I_C$} \cr
       \{x_{i-1},x_i\} & \mbox{if $x_i\in I_C^+\setminus I_C$}\cr
       \{x_i\} & \mbox{if $x_i\notin I_C\cup I_C^+$}.
     \end{array}
     \right.
$$
Then, $y_i,z_i\in V(G)-V(C)$ and
\begin{equation}\label{eqn-3-17}
E_G(\{y_i,z_i\},V(C))=\{y_ix:~x\in X_i\},~\forall ~i\in [1,m].
\end{equation}
By reversing the orientation of $C$ (if necessary), we may assume that $N(y_1)\cap V(C)\subseteq\{x_m,x_1\}$. If $m\geq 9$, then by (\ref{eqn-3-17}), we can derive that $z_1y_1x_1x_2\ldots x_8$ is an induced $P_{10}$ in $G$, a contradiction. Thus, $m\leq 8$.
As $G$ does not contain $C_8$, $m=7$.

It follows from Lemma \ref{lem-2-4} that $I_C=\emptyset$. This together with
(\ref{eqn-3-17}) implies that
\begin{equation}\label{eqn-3-18}
E_G(\{y_i,z_i\},V(C))=\{y_ix_i\},~\forall~i\in [1,7].
\end{equation}
It follows that  $N_C(y_1), N_C(y_2),\ldots,N_C(y_7)$  are distinct subsets of $V(C)$, and hence $y_1, y_2, \ldots, y_7$ are distinct vertices of $G-C$.
We claim that
\begin{equation}\label{eqn-3-19}
\{y_i,z_i\}\cap\{y_j,z_j\}=\emptyset, ~1\leq i<j\leq 7.
\end{equation}
By way of contradiction, assume that (\ref{eqn-3-19}) is false for some integers $i,j$ with $1\leq i<j\leq 7$. Then, by (\ref{eqn-3-18}), we can derive that
$\{y_i,z_i\}\cap\{y_j,z_j\}=\{z_i\}\cap\{z_j\}=\{z_i\}$. Set
$C':=x_i\overrightarrow{C}x_jy_jz_iy_ix_i$
and $C'':=x_i\overleftarrow{C}x_jy_jz_iy_ix_i$. Then, $C'$ and $C''$ are cycles  in $G$ of length at least $5$ such that
\begin{equation}\label{eqn-3-20}
\ell(C')+\ell(C'')=\ell(C)+8=15.
\end{equation}
By symmetry, we may assume that $\ell(C')\leq\ell(C'')$. This together with (\ref{eqn-3-20}) and $\ell(C'')\neq 8$ implies that $\ell(C')\leq 6<\ell(C)$, contrary to Claim \ref{claim-3-1}. Hence, (\ref{eqn-3-19}) is true. By an argument similar to that in the proof of (\ref{eqn-3-19}), we can derive that
\begin{equation}\label{eqn-3-21}
E_G(\{y_i,z_i\},\{y_j,z_j\})\subseteq\{z_iz_j\}, ~1\leq i<j\leq 7.
\end{equation}
If $z_1z_3\notin E(G)$, then by (\ref{eqn-3-18}), (\ref{eqn-3-19}) and (\ref{eqn-3-21}), we can derive that
$$
z_3y_3x_3x_4x_5x_6x_7x_1y_1z_1z_3
$$
is an induced $P_{10}$ in $G$, a contradiction. Hence, $z_1z_3\in E(G)$. Similarly, we have $z_3z_5, z_5z_7\in E(G)$. This together with (\ref{eqn-3-19}) implies that $x_1y_1z_1z_3z_5z_7y_7x_7x_1$ is a $C_8$ in $G$, a contradiction. Hence, Claim \ref{claim-3-2} is true. \qed

\begin{claim} \label{claim-3-3}
$\ell(C)\neq 6$.
\end{claim}
\noindent{\bf Proof}\ By way of contradiction, assume that $\ell(C)=6$. By Lemma \ref{lem-2-4}, we have $|I_C|\leq 1$. By permuting the indices of $x_i$ (if necessary), we may assume that $I_C\subseteq\{x_1\}$. By Lemma \ref{lem-2-5}, for each $i\in [1,6]$, $G$ has a good path, say  $y_iz_i$,  for $(C,X_i)$ with order $2$, where
$$
X_i:=\left\{\begin{array}{ll}
       \{x_i,x_{i+1}\} & \mbox{if $x_i\in I_C$} \cr
       \{x_{i-1},x_i\} & \mbox{if $x_i\in I_C^+\setminus I_C$}\cr
       \{x_i\} & \mbox{if $x_i\notin I_C\cup I_C^+$}.
     \end{array}
     \right.
$$
Then, $y_i,z_i\notin V(C)$,
\begin{equation}\label{eqn-3-22}
E_G(\{y_1,z_1\},V(C))=\left\{
  \begin{array}{ll}
   \{y_1x_1, y_1x_2\}, & \hbox{if $I_C=\{x_1\}$}\cr
   \{y_1x_1\}, & \hbox{if $I_C=\emptyset$}\cr
\end{array}
\right.
\end{equation}
and
\begin{equation}\label{eqn-3-23}
E_G(\{y_i,z_i\},V(C))=\{y_ix_{i}\}, ~i\in [3,6].
\end{equation}
It follows from (\ref{eqn-3-22}) that $N(z_1)\cap V(C)=\emptyset$.
As $\delta(G)\geq 3$, $z_1$ has two neighbors, say $u_1,u_1'$, in $V(G)-(V(C)\cup\{y_1\})$. We claim that
\begin{equation}\label{eqn-3-24}
\min~\{|N(u_1)\cap \{y_1,x_4\}|,~|N(u_1')\cap \{y_1,x_4\}|\}=0.
\end{equation}
By way of contradiction, assume  (\ref{eqn-3-24}) is false. Then, both $u_1$ and $u_1'$ has a neighbor in $\{y_1,x_4\}$. By applying Lemma \ref{lem-2-3} with $(u,v,v',A):=(z_1,u_1,u_1',\{y_1,x_4\})$, we can derive that $E_G(\{u_1,u_1'\},\{y_1,x_4\})$ contains two independent edges.
By renaming $u_1$ and $u_1'$ (if necessary), we may assume that $y_1\in N(u_1)$ and $x_4\in N(u_1')$. Then, $x_1y_1u_1z_1u_1'x_4x_5x_6x_1$ is a $C_8$ in $G$, a contradiction. Hence, (\ref{eqn-3-24}) is true.

By symmetry, we may assume that
\begin{equation}\label{eqn-3-25}
N(u_1)\cap \{y_1,x_4\}=\emptyset.
\end{equation}
If $N(u_1)\cap (V(C)\cup \{y_1\})\neq\emptyset$, then by (\ref{eqn-3-25}), we see that
$x_i\in N(u_1)$ holds for some $i\in\{1,2,3,5,6\}$. If $i\in\{1,3,5\}$, then
$$
C':=\left\{
  \begin{array}{ll}
   x_1y_1z_1u_1x_1, & \hbox{if $i=1$}\cr
   x_1y_1z_1u_1x_3x_4x_5x_6x_1, & \hbox{if $i=3$}\cr
   x_1y_1z_1u_1x_5x_4x_3x_2x_1 & \hbox{if $i=5$}
\end{array}
\right.
$$
is a $C_4$ or $C_8$ in $G$, a contradiction. Hence, $i\notin\{1,3,5\}$, which in turn means $i\in\{2,6\}$. It follows that $x_1y_1z_1u_1x_ix_1$ is a $C_5$ in $G$, contrary to Claim \ref{claim-3-1}.  Hence,
\begin{equation}\label{eqn-3-26}
N(u_1)\cap (V(C)\cup \{y_1\})=\emptyset.
\end{equation}

Let $i\in\{1,3,5\}$ and define $(x_j,y_j,z_j):=(x_{j-6},y_{j-6},z_{j-6})$, where $j\in [7,12]$. It follows from (\ref{eqn-3-22}) and (\ref{eqn-3-23}) that $N_C(z_i)=\emptyset$, and
$$
N_C(y_i)=\left\{
  \begin{array}{ll}
   \{x_1, x_2\}, & \hbox{if $i=1$ and $I_C=\{x_1\}$}\cr
   \{x_i\}, & \hbox{otherwise}.
\end{array}
\right.
$$
Thus, $\{y_i,z_i\}\cap \{y_{i+2},z_{i+2}\}\subseteq\{z_i\}\cap \{z_{i+2}\}$. If
$z_i=z_{i+2}$, then $x_iy_iz_iy_{i+2}x_{i+2}\overrightarrow{C}x_i$ is a $C_8$ in $G$, a contradiction. Hence, $z_i\neq z_{i+2}$. It follows that
\begin{equation}\label{eqn-3-27}
\{y_i,z_i\}\cap\{y_{i+2},z_{i+2}\}=\emptyset,~i=1,3,5.
\end{equation}
By an argument similar to that in the proof of $z_i\neq z_{i+2}$, we can derive that $z_iy_{i+2},y_iz_{i+2}\notin E(G)$. If $y_iy_{i+2}\in E(G)$, then $x_iy_iy_{i+2}x_{i+2}x_{i+1}x_i$ is a $C_5$ in $G$, contrary to Claim \ref{claim-3-1}. Hence, $y_iy_{i+2}\notin E(G)$. Therefore,
\begin{equation}\label{eqn-3-28}
E_G(\{y_i,z_i\},\{y_{i+2},z_{i+2}\})\subseteq\{z_iz_{i+2}\},~i=1,3,5.
\end{equation}
We claim that
\begin{equation}\label{eqn-3-29}
|E(G)\cap\{z_1z_3, z_3z_5, z_5z_1\}|\leq 1.
\end{equation}
Otherwise, there exists $i\in \{1,3,5\}$ such that $\{z_iz_{i+2}, z_{i+2}z_{i+4}\}\subseteq E(G)$. It follows that $x_iy_iz_iz_{i+2}z_{i+4}y_{i+4}x_{i+4}
x_{i+5}x_i$ is a $C_8$ in $G$, a contradiction. Hence, (\ref{eqn-3-29}) is true.
\begin{figure}\label{fig-6}
\psfrag{A}{$x_i$}
\psfrag{B}{$x_{i+1}$}
\psfrag{C}{$x_{i+2}$}
\psfrag{D}{$x_{i+3}$}
\psfrag{E}{$x_{i+4}$}
\psfrag{F}{$x_{i+5}$}
\psfrag{G}{$y_{i}$}
\psfrag{H}{$z_i$}
\psfrag{I}{$y_{i+2}$}
\psfrag{J}{$z_{i+2}$}
\psfrag{K}{$y_{i+4}$}
\psfrag{L}{$z_{i+4}$}
\psfrag{M}{$I_C\neq\{x_{i+2}\}$}
\psfrag{N}{$I_C=\{x_{i+2}\}$}
\begin{center}
     \includegraphics[width=130mm]{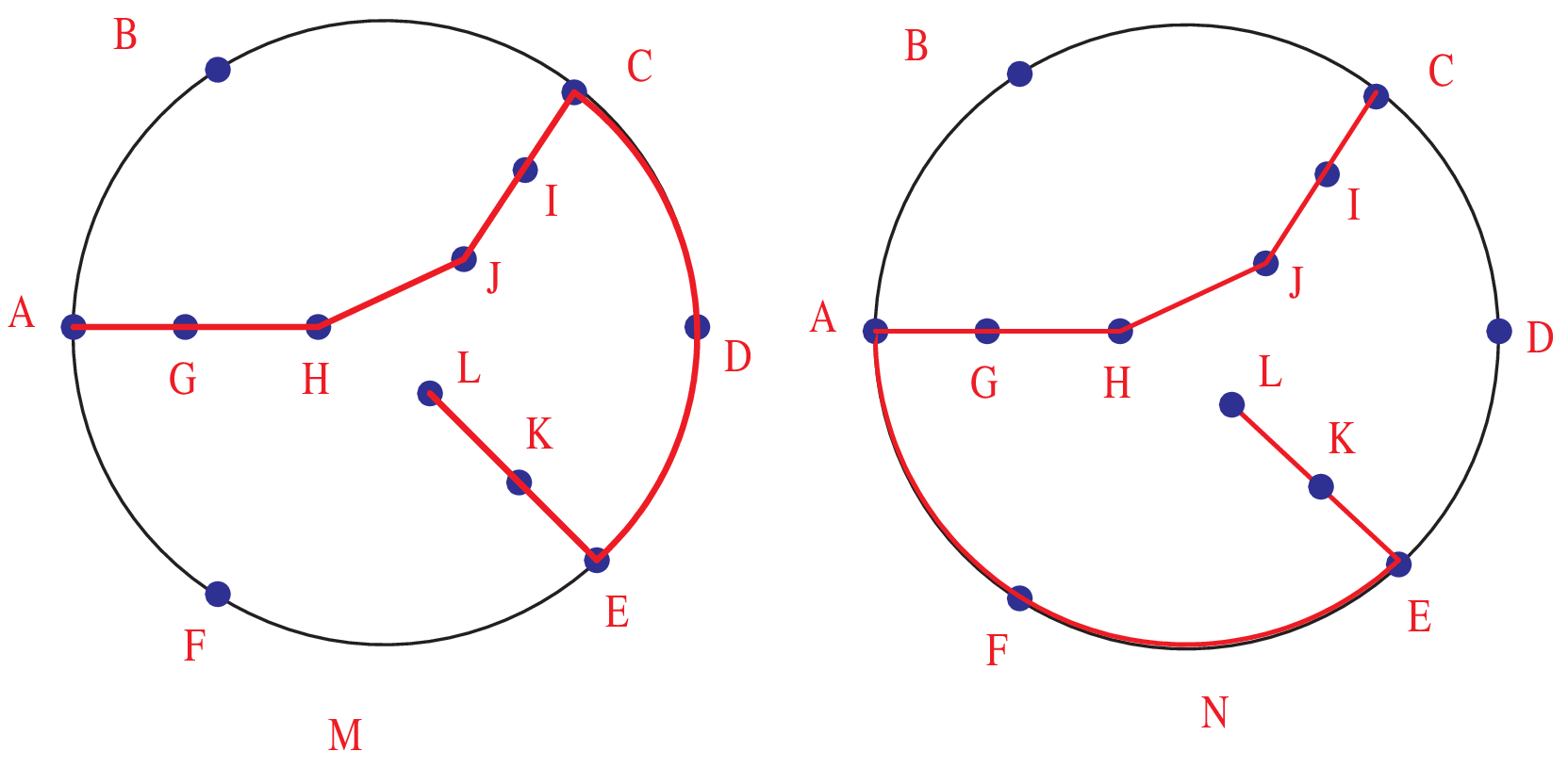}\\
   \caption{Illustration of the nonexistence of the edge $z_iz_{i+2}$.}
   \end{center}
\end{figure}
Now, we claim that
\begin{equation}\label{eqn-3-30}
E_G(\{y_i,z_i\},\{y_{i+2},z_{i+2}\})=\emptyset,~i=1,3,5.
\end{equation}
By way of contradiction, assume that  (\ref{eqn-3-30}) is false for some integer $i\in\{1,3,5\}$. Then, by  (\ref{eqn-3-28}) and  (\ref{eqn-3-29}), we may assume that
$$
E_G(\{y_j,z_j\},\{y_{j+2},z_{j+2}\})=\left\{
  \begin{array}{ll}
   \{z_iz_{i+2}\}, & \hbox{if $j=i$}\cr
   \emptyset, & \hbox{if $j=i+2,i+4$}.
\end{array}
\right.
$$
This together with (\ref{eqn-3-22}) and  (\ref{eqn-3-23}) implies that
$$
R:=\left\{
  \begin{array}{ll}
   x_iy_iz_iz_{i+2}y_{i+2}x_{i+2}x_{i+3}x_{i+4}y_{i+4}z_{i+4}, & \hbox{if $I_C\neq\{x_{i+2}\}$}\cr
  z_{i+4}y_{i+4}x_{i+4} x_{i+5}x_iy_iz_iz_{i+2}y_{i+2}x_{i+2}, & \hbox{if $I_C=\{x_{i+2}\}$}.
\end{array}
\right.
$$
is an induced $P_{10}$ in $G$ (see Figure 6), a contradiction.
Hence, (\ref{eqn-3-30}) is true. 

Recall that $u_1\in N(z_1)-(V(C)\cup\{y_1\})$. By (\ref{eqn-3-30}), we can derive  that $y_1,z_1,u_1,y_3$, $z_3,y_5,z_5$ are distinct vertices of $G-C$. If $u_1\notin N(y_3)\cup N(z_3)$, then by  (\ref{eqn-3-22}), (\ref{eqn-3-23}), (\ref{eqn-3-26}) and  (\ref{eqn-3-30}), we can derive  that $z_3y_3x_3x_4x_5x_6x_1y_1z_1u_1$ is an induced $P_{10}$ in $G$, a contradiction. Hence,
\begin{equation}\label{eqn-3-31}
u_1\in N(y_3)\cup N(z_3).
\end{equation}
If $u_1\in N(z_3)$, then $x_1y_1z_1u_1z_3y_3x_3x_2x_1$ is a $C_8$ in $G$, a contradiction. Hence, $u_1\notin N(z_3)$.  By combining this with (\ref{eqn-3-31}), we get
$u_1\in N(y_3)$.  This together with (\ref{eqn-3-23}) and  (\ref{eqn-3-26}) implies that $y_3u_1$ is a good path for $(C,X_3)$. Recall that  $y_1z_1$ is a good path for $(C,X_1)$. By an argument similar to that in the proof of  (\ref{eqn-3-30}), we can derive that
$E_G(\{y_1,z_1\},\{y_3,u_1\})=\emptyset$,   contrary to $u_1\in N(z_1)$. Hence, Claim \ref{claim-3-3} is true.\qed

\vspace{3mm}
It follows from Claims \ref{claim-3-2} and \ref{claim-3-3} that $\ell(C)=5$.
In order to complete the proof of Theorem \ref{thm1.1}, we show some  claims about the structure of $G-C$, where the indices of $x_i$ and $A_i$ are taken modulo $5$.

\begin{claim}\label{cla-3-4}
Let $i,j$ be two distinct integers in $[1,5]$ and let $P$ be a $(u,v)$-path in $G-V(C)$ with $u\in A_i$ and $v\in A_j$. Then

\emph{(i)} if $j\equiv i\pm 1~(mod~5)$, then $\ell(P)\neq 1,2,5$;  and

\emph{(ii)} if $j\equiv i\pm 2~(mod~5)$, then $\ell(P)\neq 0,3,4$.
\end{claim}

\noindent{\bf Proof}\  By way of contradiction, assume that Claim  \ref{cla-3-4} is false. By reversing the orientation of $C$ (if necessary), we may assume that $j= i+1,i+2$. Then $\ell(x_i\overrightarrow{C}x_j)=j-i$ and $\ell(x_i\overleftarrow{C}x_j)=5-(j-i)$. It follows that
$$
C':=\left\{\begin{array}{ll}
x_i\overrightarrow{C}x_jv{P}ux_i,&\makebox{if $(j-i,\ell(P))=(1,1),(1,5),(2,0),(2,4)$}\cr
x_i\overleftarrow{C}x_jv{P}ux_i,&\makebox{if $(j-i,\ell(P))=(1,2),(2,3)$}
\end{array}\right.
$$
is a $C_4$ or $C_8$ in $G$, a contradiction. Hence, Claim  \ref{cla-3-4} is true.\qed

\begin{claim}\label{cla-3-5}
Let $x_i\in I_C$ and let $P:=u_1\ldots u_s$, $Q:=v_1\ldots v_t$ be two paths in $G-C$ such that $u_1\in A_i\cap A_{i+1}$ and $v_1\in A_{i+2}\cup A_{i+4}$.

\emph{(i)} If $s+t\leq 7$, then $V(P)\cap V(Q)=\emptyset$; and

\emph{(ii)} if $s+t\leq 6$, then $E_G(V(P),V(Q))=\emptyset$.
\end{claim}

\noindent{\bf Proof}\ By way of contradiction, assume that Claim \ref{cla-3-5} is false. Then $G[V(P)\cup V(Q)]$ contains a $(u_1,v_1)$-path $R$ of length at most five with $u_1\in A_i\cap A_{i+1}$ and $v_1\in A_{i+2}\cup A_{i+4}$ such that $V(R)\subseteq V(G)-V(C)$.  Note that $A_{i+4}=A_{i-1}$. We may assume, by symmetry, that $v_1\in A_{i+2}$.
\begin{itemize}
  \item If $\ell(R)\in\{1,2,5\}$, then by Claim \ref{cla-3-4} (i), we can derive that
  $u_1\notin A_{i+1}$;
  \item If $\ell(R)\in\{0,3,4\}$, then by Claim \ref{cla-3-4} (ii), we have $u_1\notin A_i$.
\end{itemize}
In both cases, we get a contradiction. Hence, Claim \ref{cla-3-5} is true.\qed

\vspace{4mm}

An \emph{$x_i$-path} in $G$ is a path starting from $x_i$. If $P$ is an $x_i$-path in $G$ such that $V(P)\cap V(C)=\{x_i\}$, then we call $P$ a \emph{pendent $x_i$-path for $C$}.

\begin{claim}\label{cla-3-6}
Let $i\in [1,5]$ and let $P:=x_iuvw$ be a pendent $x_i$-path for $C$.
Then $x_{i-1},x_i,x_{i+1}\notin N(w)$.
\end{claim}
\noindent{\bf Proof}\ By way of contradiction, assume that Claim \ref{cla-3-6} is false. Then, $w$ has a neighbor, say $z$, in $\{x_{i-1},x_i,x_{i+1}\}$. Set
$$
C':=\left\{\begin{array}{ll}
x_iuvwx_{i-1}\overleftarrow{C}x_i,&\makebox{if $z=x_{i-1}$}\cr
x_iuvwx_i,&\makebox{if $z=x_{i}$}\cr
x_iuvwx_{i+1}\overrightarrow{C}x_i,&\makebox{if $z=x_{i+1}$}.
\end{array}\right.
$$
Then, $C'$ is a $C_4$ or $C_8$ in $G$, a contradiction.  Hence, Claim \ref{cla-3-6} is true.\qed

\begin{defi}\label{def-3-1}
Let $w_1\ldots w_{\ell}$ be an  induced path of $G-C$ with length $2$ or $3$. If there exists an integer $i\in [1,5]$ such that
\begin{itemize}
  \item $N_C(w_1)=\{x_i,x_{i+1}\}$, $N_C(w_2)\subseteq \{x_{i+3}\}$, $N_C(w_3)=\emptyset$, and
  \item $N_C(w_4)\subseteq\{x_i,x_{i+1}\}$ if $\ell=4$,
\end{itemize}
then we call $x_iw_1\ldots w_{\ell}$ a \emph{near-good $(C,x_i)$-path} of length $\ell$.
\end{defi}

\begin{claim}\label{cla-3-7}
For each  $x_i\in I_C$, there exists a near-good $(C,x_i)$-path.
\end{claim}
\noindent{\bf Proof}\ Let $x_i\in I_C$. Then, $A_i\cap A_{i+1}\neq\emptyset$.
We will find a near-good $(C,x_i)$-path $x_iy_iz_iu_iv_i$ step by step.

First, we let $y_i\in A_i\cap A_{i+1}$. By the proof of Lemma \ref{lem-2-5}, we see that $y_i$ is a good $(C,\{x_i,x_{i+1}\})$-path. Thus,
\begin{equation}\label{eqn-3-32}
    N_C(y_i)=\{x_i,x_{i+1}\}.
\end{equation}

Note that $\delta(G)\geq 3$. By (\ref{eqn-3-32}),  we have $N(y_i)-V(C)\neq\emptyset$.
Choose $z_i\in N(y_i)-V(C)$ such that $|N_C(z_i)|$ is as small as possible. We claim that
\begin{equation}\label{eqn-3-33}
    N_C(z_i)\subseteq\{x_{i+3}\}.
\end{equation}
Otherwise, $z_ix_j\in E(G)$ holds for some $j\in\{i-1,i,i+1,i+2\}$. By symmetry, we may assume that $j\in\{i+1,i+2\}$. Then, $y_ix_{j-1}\in E(G)$. It follows that $x_{j-1}y_iz_ix_jx_{j-1}$ is a $C_4$ in $G$, a contradiction. Hence, (\ref{eqn-3-33}) is true.

It follows from (\ref{eqn-3-33}) that
$N(z_i)-(V(C)\cup\{y_i\})= N(z_i)-\{x_{i+3},y_i\}$.
As $\delta(G)\geq 3$, $N(z_i)-(V(C)\cup\{y_i\})\neq\emptyset$.
Among all vertices of $N(z_i)-(V(C)\cup\{y_i\})$, choose one, say $u_i$, such that $|N(u_i)\cap \{x_{i+3},y_i\}|$ achieves the minimum. We claim that
\begin{equation}\label{eqn-3-34}
  N(u_i)\cap \{x_{i+3},y_i\}=\emptyset.
\end{equation}
By way of contradiction, assume that (\ref{eqn-3-34}) is false. We consider two cases.
\begin{itemize}
  \item $N(z_i)-(V(C)\cup\{y_i\})=\{u_i\}$. As $\delta(G)\geq 3$,  we have $N_C(z_i)\neq\emptyset$. This together with (\ref{eqn-3-33}) implies that
      $N_C(z_i)=\{x_{i+3}\}$.  As (\ref{eqn-3-34}) is false, $N(u_i)\cap\{x_{i+3},y_i\}\neq\emptyset$. If $N(u_i)\supseteq \{x_{i+3},y_i\}$, then $y_iz_ix_{i+3}u_iy_i$ is a $C_4$ in $G$ (see Figure 7), a contradiction. Hence,
      $$
      |N(u_i)\cap \{x_{i+3},y_i\}|=1.
      $$
      Set $C':=x_iy_iz_ix_{i+3}x_{i+4}x_i$ and $C'':=x_{i+1}x_{i+2}x_{i+3}z_iy_ix_{i+1}$.
      Then, both $C'$ and $C''$ are $5$-holes in $G$. Note that
$$
I_{C'}=
\left\{\begin{array}{ll}
(I_C\cap\{x_{i+3},x_{i+4}\})\cup \{x_i,z_i\}, &\mbox{if $x_{i+3}\in N(u_i)$}\cr
(I_C\cap\{x_{i+3},x_{i+4}\})\cup \{x_i,y_i\}, &\mbox{if $y_i\in N(u_i)$}
\end{array}\right.
$$
and
$$
I_{C''}=
\left\{\begin{array}{ll}
(I_C\cap\{x_{i+1},x_{i+2}\})\cup \{y_i,x_{i+3}\}, &\mbox{if $x_{i+3}\in N(u_i)$}\cr
(I_C\cap\{x_{i+1},x_{i+2}\})\cup \{y_i,z_i\}, &\mbox{if $y_i\in N(u_i)$}.
\end{array}\right.
$$
Hence, $|I_{C'}|+|I_{C''}|=|I_C|+3$. On the other hand,
by the choice of $C$, we have $|I_C|\geq |I_{C'}|$, and hence $|I_{C''}|\geq 3$,
contrary to Lemma \ref{lem-2-4}.

\begin{figure}\label{fig-7}
\psfrag{A}{$x_i$}
\psfrag{B}{$x_{i+1}$}
\psfrag{C}{$x_{i+2}$}
\psfrag{D}{$x_{i+3}$}
\psfrag{E}{$x_{i+4}$}
\psfrag{F}{$y_i$}
\psfrag{G}{$z_i$}
\psfrag{H}{$u_i$}
      \begin{center}
     \includegraphics[width=60mm]{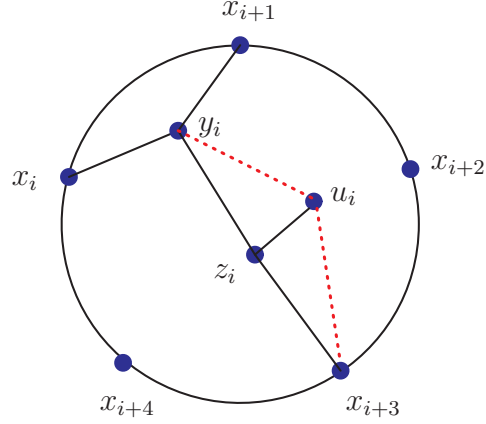}\\
   \caption{Exactly one of $x_{i+3}$ and $y_i$ is a neighbor of $u_i$.}
   \end{center}
\end{figure}

\item $N(z_i)-(V(C)\cup\{y_i\})\neq\{u_i\}$. Let $u_i'$  be a vertex of $N(z_i)-(V(C)\cup\{y_i\})$ with $u_i'\neq u_i$. As (\ref{eqn-3-34}) is false, we have
$$
   \min~\{|N(u_i)\cap \{x_{i+3},y_i\}|, |N(u_i')\cap \{x_{i+3},y_i\}|\}>0.
$$
By applying Lemma \ref{lem-2-3} with $(u,v,v',A):=(z_i,u_i,u_i',\{x_{i+3},y_i\})$, we see that $E_G(\{u_i,u_i'\},\{x_{i+3},y_i\})$ contains two independent edges. By symmetry, we may assume that $u_ix_{i+3}, u_i'y_i\in E(G)$. Then, $x_iy_iu_i'z_iu_ix_{i+3}x_{i+2}x_{i+1}x_i$ is a $C_8$ in $G$, a contradiction.
\end{itemize}
In both cases, we get a contradiction. Hence, (\ref{eqn-3-34}) is true.

Note that both $x_iy_iz_iu_i$ and $x_{i+1}y_iz_iu_i$ are pendent paths for $C$ with length three. By Claim \ref{cla-3-6}, we can derive that $x_{i-1},x_i,x_{i+1},x_{i+2}\notin N(u_i)$. This together with (\ref{eqn-3-34}) implies that
\begin{equation}\label{eqn-3-35}
    N(u_i)\cap (V(C)\cup\{y_i\})=\emptyset.
\end{equation}

It follows from (\ref{eqn-3-35}) that $|N(u_i)-(V(C)\cup\{y_i,z_i\})|=|N(u_i)-\{z_i\})|\geq 2$.
Let $v_i, v_i'$ be two vertices of $N(u_i)-(V(C)\cup\{y_i,z_i\})$.
If $v_i,v_i'\in N(z_i)$, then $v_iz_iv_i'u_iv_i$ is a $C_4$ in $G$, a contradiction.
Hence, $\{v_i,v_i'\}\not\subseteq N(z_i)$. By symmetry, we may assume that
$v_i\notin N(z_i)$. We claim that
\begin{equation}\label{eqn-3-36}
    N(v_i)\cap (V(C)\cup\{y_i,z_i\})\subseteq\{x_i,x_{i+1}\}.
\end{equation}
Otherwise, $wv_i\in E(G)$ holds for some vertex $w\in \{x_{i+2},x_{i+3},x_{i+4},y_i\}$.
Set
$$
C^*:=\left\{\begin{array}{ll}
x_iy_iz_iu_iv_ix_{i+2}x_{i+3}x_{i+4}x_i,&\makebox{if $w=x_{i+2}$}\cr
x_ix_{i+1}y_iz_iu_iv_ix_{i+3}x_{i+4}x_i,&\makebox{if $w=x_{i+3}$}\cr
x_{i+1}y_iz_iu_iv_ix_{i+4}x_{i+3}x_{i+2}x_{i+1},&\makebox{if $w=x_{i+4}$}\cr
y_iz_iu_iv_iy_i,&\makebox{if $w=y_i$}.
\end{array}\right.
$$
Then, $C^*$ is a $C_4$ or $C_8$ in $G$, a contradiction.  Hence,  (\ref{eqn-3-36}) is true.

It follows from (\ref{eqn-3-32}),(\ref{eqn-3-33}),(\ref{eqn-3-35}) and (\ref{eqn-3-36}) that $x_iy_iz_iu_iv_i$ is a near-good $(C,x_i)$-path. Hence, Claim \ref{cla-3-7} is true. \qed

\begin{claim}\label{cla-3-8}
Let $i$ be an integer in $[1,5]$ such that $x_i\notin I_C\cup I_C^+$, then
there exists a good $(C,x_i)$-path of length three in $G$.
\end{claim}
\noindent{\bf Proof}\  It follows from Lemma \ref{lem-2-5} that $G$ contains a good $(C,x_i)$-path, say $x_iy_i$, of length one.
As $x_i\notin I_C\cup I_C^+$, $E_G(\{y_i\}, V(C))=\{y_ix_i\}$. Thus,
\begin{equation}\label{eqn-3-37}
    N_C(y_i)=\{x_i\}.
\end{equation}
As $\delta(G)\geq 3$, $y_i$ has at least two neighbors, say $z_i$ and $z_i'$, in $V(G)-V(C)$.  We claim that
\begin{equation}\label{eqn-3-38}
    \min~\{|N(z_i)\cap \{x_i, x_{i+2},x_{i+3}\}|, |N(z_i')\cap \{x_i, x_{i+2},x_{i+3}\}|\}=0.
\end{equation}
By way of contradiction, assume that (\ref{eqn-3-38}) is false.
By applying Lemma \ref{lem-2-3} with $(u,v,v',A):=(y_i,z_i,z_i',\{x_i, x_{i+2},x_{i+3}\})$, we see that there exists two independent edges, say $e_1$ and $e_2$, in $E_G(\{z_i,z_i'\},\{x_i, x_{i+2},x_{i+3}\})$ such that $e_1$ has an end vertex in $\{x_{i+2},x_{i+3}\}$. Based on $(C,x_iy_i)$, $z_i$ and $z_i'$ are symmetrical, while
$x_{i+2}$ and $x_{i+3}$ are symmetrical (see Figure 8).
\begin{figure}\label{fig-8}
\psfrag{A}{$x_i$}
\psfrag{B}{$x_{i+1}$}
\psfrag{C}{$x_{i+2}$}
\psfrag{D}{$x_{i+3}$}
\psfrag{E}{$x_{i+4}$}
\psfrag{F}{$y_i$}
\psfrag{G}{$z_i$}
\psfrag{H}{$z_i'$}
\psfrag{I}{$e_1$}
\psfrag{J}{$e_2$}
      \begin{center}
     \includegraphics[width=60mm]{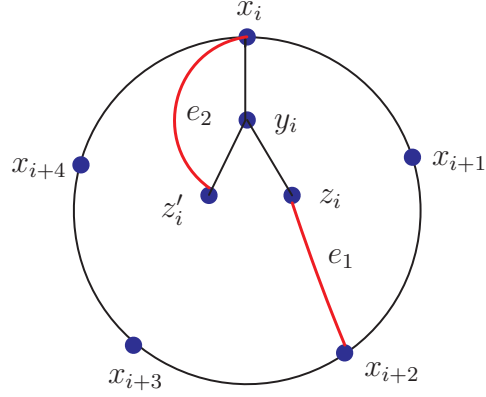}\\
   \caption{The case $\{e_1,e_2\}=\{z_ix_{i+2},z_i'x_i\}$.}
   \end{center}
\end{figure}
Thus, we may assume that $e_1=z_ix_{i+2}$, which in turn means $e_2\in\{z_i'x_i,z_i'x_{i+3}\}$.

\begin{itemize}
  \item If $e_2=z_i'x_{i+3}$, then
$x_{i+2}z_iy_iz_i'x_{i+3}x_{i+4}x_ix_{i+1}x_{i+2}$ is a $C_8$ in $G$, a contradiction.
\item If $e_2=z_i'x_i$, set $C':=x_ix_{i+1}x_{i+2}z_iy_ix_i$ and $C'':=x_iy_iz_ix_{i+2}x_{i+3}x_{i+4}x_i$. Then, $C'$ is a $5$-hole in $G$ with $I_{C'}\supseteq(I_C\cap\{x_{i+1}\})\cup\{y_i\}$. By the choice of $C$, we have $|I_{C}|\geq|I_{C'}|$, and hence $I_C\setminus\{x_{i+1}\}\neq\emptyset$. This together with $x_i\notin I_C\cup I_C^+$ implies that $I_C\cap\{x_{i+2},x_{i+3}\}\neq\emptyset$. Thus, $x_j\in I_C$ holds for some $j\in\{i+2,i+3\}$. It follows that $C''$ is a $C_6$ in $G$ such that
$I_{C''}\supseteq\{x_i,x_j\}$. On the other hand, by Lemma \ref{lem-2-4}, we have $|I_{C''}|\leq 7-6$, a contradiction.
\end{itemize}
Therefore, (\ref{eqn-3-38}) is true.

By symmetry, we may assume that $N(z_i)\cap \{x_i, x_{i+2},x_{i+3}\}=\emptyset$.
If $N_C(z_i)\neq\emptyset$, then $z_i$ has a neighbor, say $w$, in $\{x_{i+1},x_{i+4}\}$. It follows that $x_iy_iz_iwx_i$ is a $C_4$ in $G$, a contradiction.  Hence,
\begin{equation}\label{eqn-3-39}
N_C(z_i)=\emptyset.
\end{equation}

It follows from (\ref{eqn-3-39}) that $|N(z_i)-(V(C)\cup\{y_i\})|=d(z_i)-1\geq 2$. Let $u_i$ and $u_i'$ be two neighbors of $z_i$ in $V(G)-(V(C)\cup\{y_i\})$. We claim that
\begin{equation}\label{eqn-3-40}
    \min~\{|N(u_i)\cap \{y_i, x_{i+2},x_{i+3}\}|, |N(u_i')\cap \{y_i, x_{i+2},x_{i+3}\}|\}=0.
\end{equation}
By way of contradiction, assume that (\ref{eqn-3-40}) is false.
By applying Lemma \ref{lem-2-3} with $(u,v,v',A):=(z_i,u_i,u_i',\{y_i, x_{i+2},x_{i+3}\})$, we see that there exists two independent edges, say $f_1$ and $f_2$, in $E_G(\{u_i,u_i'\},\{y_i, x_{i+2},x_{i+3}\})$ such that $f_1$ has an end vertex in $\{x_{i+2},x_{i+3}\})$. Based on $(C,x_iy_iz_i)$, $u_i$ and $u_i'$ are symmetrical, while
$x_{i+2}$ and $x_{i+3}$ are symmetrical (see Figure 9).
\begin{figure}\label{fig-9}
\psfrag{A}{$x_i$}
\psfrag{B}{$x_{i+1}$}
\psfrag{C}{$x_{i+2}$}
\psfrag{D}{$x_{i+3}$}
\psfrag{E}{$x_{i+4}$}
\psfrag{F}{$y_i$}
\psfrag{G}{$z_i$}
\psfrag{H}{$u_i$}
\psfrag{I}{$u_i'$}
\psfrag{J}{$f_1$}
\psfrag{K}{$f_2$}
      \begin{center}
     \includegraphics[width=60mm]{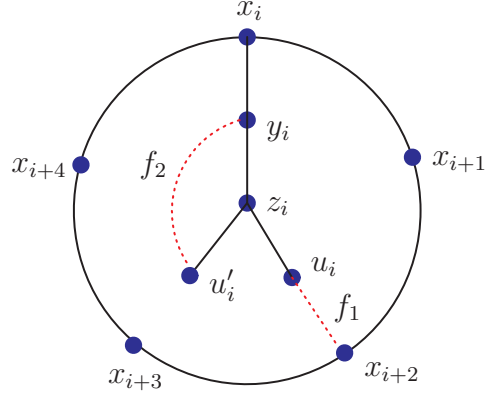}\\
   \caption{The case $(f_1,f_2)=(u_ix_{i+2},u_i'y_i)$.}
   \end{center}
\end{figure}
Thus, we may assume that $f_1=u_ix_{i+2}$, which in turn means $f_2\in\{u_i'y_i,u_i'x_{i+3}\}$. It follows that
$$
D':=\left\{\begin{array}{ll}
  x_iy_iu_i'z_iu_ix_{i+2}x_{i+3}x_{i+4}x_i, &\mbox{if $f_2=u_i'y_i$}\cr
  x_{i+2}u_iz_iu_i'x_{i+3}x_{i+4}x_ix_{i+1}x_{i+2}, &\mbox{if $f_2=u_i'x_{i+3}$}
\end{array}\right.
$$
is a $C_8$ in $G$, a contradiction. Hence, (\ref{eqn-3-40}) is true.

By symmetry, we may assume that $N(u_i)\cap \{y_i, x_{i+2},x_{i+3}\}=\emptyset$.
If $N(u_i)\cap (V(C)\cup\{y_i\})\neq\emptyset$, then $u_i$ has a neighbor, say $w$, in $\{x_i, x_{i+1},x_{i+4}\}$. It follows that
$$
D'':=\left\{\begin{array}{ll}
  x_iy_iz_iu_ix_i, &\mbox{if $w=x_i$}\cr
  x_iy_iz_iu_ix_{i+1}x_{i+2}x_{i+3}x_{i+4}x_i, &\mbox{if $w=x_{i+1}$}\cr
  x_iy_iz_iu_ix_{i+4}x_{i+3}x_{i+2}x_{i+1}x_i, &\mbox{if $w=x_{i+4}$}
\end{array}\right.
$$
is a $C_4$ or $C_8$ in $G$, a contradiction.  Hence,
\begin{equation}\label{eqn-3-41}
N(u_i)\cap (V(C)\cup\{y_i\})=\emptyset.
\end{equation}
By (\ref{eqn-3-37}),(\ref{eqn-3-39}) and (\ref{eqn-3-41}), we see that
$x_iy_iz_iu_i$ is a good $(C,x_i)$-path of length three in $G$.
Hence, Claim \ref{cla-3-8} is true. \qed

\begin{claim}\label{cla-3-9}
For each  $x_i\in I_C$, there exists a near-good $(C,x_i)$-path $x_iy_iz_iu_iv_i$ such that
$$
E_G(\{y_i,z_i,u_i,v_i\}, V(C))=\{y_ix_i,y_ix_{i+1},z_ix_{i+3}\}.
$$
\end{claim}
\noindent{\bf Proof}\ Let $x_i\in I_C$. By Claim \ref{cla-3-7}, $G$ contains  a near-good $(C,x_i)$-path, say $x_iy_iz_iu_iv_i$, of length $4$.
It follows from  Lemma \ref{lem-2-4} that $|I_C|\leq 7-5$, and hence $|I_C\cap\{x_{i+1},x_{i+2}\}|+|I_C\cap\{x_{i-1},x_{i-2}\}|\leq 1$.
By symmetry, we may assume that $I_C\cap\{x_{i+1},x_{i+2}\}=\emptyset$. Then, $x_{i+2}\notin I_C\cup I_C^+$. By Claim \ref{cla-3-8}, $G$ contains a good $(C,x_{i+2})$-path, say $x_{i+2}y_{i+2}z_{i+2}u_{i+2}$, of length $3$. It follows from Definitions \ref{def-2-2} and \ref{def-3-1} that both $y_iz_iu_iv_i$ and $y_{i+2}z_{i+2}u_{i+2}$ are induced paths of $G-C$ such that
\begin{equation}\label{eqn-3-42}
\{y_ix_i,y_ix_{i+1}\}\subseteq E_G(\{y_i,z_i,u_i,v_i\}, V(C))\subseteq \{y_ix_i,y_ix_{i+1},z_ix_{i+3},v_ix_i,v_ix_{i+1}\}
\end{equation}
and
\begin{equation}\label{eqn-3-43}
E_G(\{y_{i+2},z_{i+2},u_{i+2}\}, V(C))= \{y_{i+2}x_{i+2}\}.
\end{equation}
Note that $y_i\in A_i\cap A_{i+1}$ and $y_{i+2}\in A_{i+2}$. By applying Claim \ref{cla-3-5} with $P:=y_iz_iu_iv_i$ and $Q:=y_{i+2}z_{i+2}u_{i+2}$, we have
\begin{equation}\label{eqn-3-44}
\{y_i,z_i,u_i,v_i\}\cap\{y_{i+2},z_{i+2},u_{i+2}\}=\emptyset
\end{equation}
and
\begin{equation}\label{eqn-3-45}
E_G(\{y_i,z_i,u_i,v_i\}, \{y_{i+2},z_{i+2},u_{i+2}\})\subseteq \{v_iu_{i+2}\}.
\end{equation}
If $z_ix_{i+3}\notin E(G)$, then by (\ref{eqn-3-42})-(\ref{eqn-3-45}), we can derive that
$$
u_{i+2}z_{i+2}y_{i+2}x_{i+2}x_{i+3}x_{i+4}x_iy_iz_iu_i
$$
is an induced $P_{10}$ in $G$ (see Figure 10), a contradiction. Hence, $z_ix_{i+3}\in E(G)$. If Claim \ref{cla-3-9} is not true, then by (\ref{eqn-3-42}), we can derive that $\{v_ix_i,v_ix_{i+1}\}\cap E(G)\neq\emptyset$. If $v_ix_i\in E(G)$, then $x_iy_ix_{i+1}x_{i+2}x_{i+3}z_iu_iv_ix_i$ is a $C_8$ in $G$;
if $v_ix_{i+1}\in E(G)$, then $x_iy_ix_{i+1}v_iu_iz_ix_{i+3}x_{i+4}x_i$ is a $C_8$ in $G$; either way gives a contradiction. Hence,  Claim \ref{cla-3-9} is  true.\qed
\begin{figure}\label{fig-10}
\psfrag{A}{$x_i$}
\psfrag{B}{$x_{i+1}$}
\psfrag{C}{$x_{i+2}$}
\psfrag{D}{$x_{i+3}$}
\psfrag{E}{$x_{i+4}$}
\psfrag{F}{$y_i$}
\psfrag{G}{$z_i$}
\psfrag{H}{$u_i$}
\psfrag{I}{$v_i$}
\psfrag{J}{$y_{i+2}$}
\psfrag{K}{$z_{i+2}$}
\psfrag{L}{$u_{i+2}$}
\psfrag{M}{Possible edges of $G$}
      \begin{center}
     \includegraphics[width=70mm]{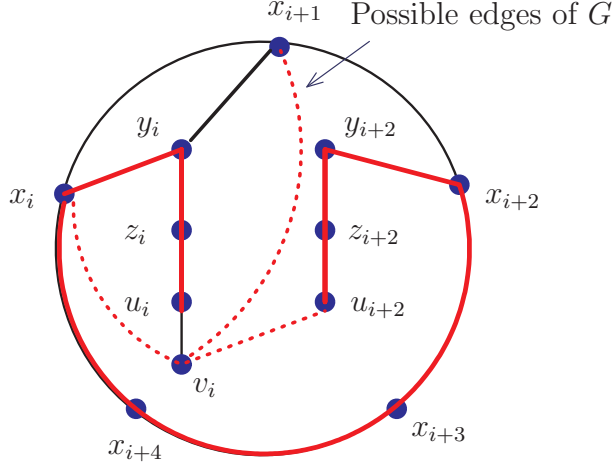}\\
   \caption{Induced $P_{10}$ in case $z_ix_{i+3}\notin E(G)$.}
   \end{center}
\end{figure}

\begin{claim}\label{cla-3-10}
$I_C=\emptyset$.
\end{claim}
\noindent{\bf Proof}\  It follows from  Lemma \ref{lem-2-4} that $|I_C|\leq 7-5$. If  Claim \ref{cla-3-10} is not true, then $1\leq |I_C|\leq 2$.
As $|V(C)|= 5$, there exists $i_0\in [1,5]$ such that $I_C\cap\{x_{i_0},x_{i_0+1},x_{i_0+2}\}=\{x_{i_0}\}$. By permuting the vertices of $C$ (if necessary), we may assume that $i_0=1$, that is $I_C\cap\{x_1,x_2,x_3\}=\{x_1\}$.
Recall that  $|I_C|\leq 2$. There are following two cases.
\begin{itemize}
  \item $I_C=\{x_1,x_k\}$, where $k\in\{4,5\}$.   For  $i=1,k$, Let $x_iy_iz_iu_iv_i$ be a near-good $(C,x_i)$-path satisfying Claim \ref{cla-3-9}. Then,
      \begin{equation}\label{eqn-3-46}
      E_G(\{y_i,z_i,u_i,v_i\}, V(C))=\{y_ix_i,y_ix_{i+1},z_ix_{i+3}\}, ~i=1,k.
      \end{equation}
      Note that $y_k\in A_k\cap A_{k+1}$ and $y_1\in A_{k+2}$. By applying Claim \ref{cla-3-5} with $P:=y_kz_ku_kv_k$ and $Q:=y_1z_1u_1$, we can derive that $\{y_k,z_k,u_k,v_k\}\cap\{y_1,z_1,u_1\}=\emptyset$.
      This together with (\ref{eqn-3-46}) implies that
 $$
 C'=\left\{\begin{array}{ll}
   x_1y_1z_1x_4x_5y_4z_4x_2x_1, &\makebox{if $k=4$}\cr
   x_2y_1z_1x_4x_3z_5y_5x_1x_2, &\makebox{if $k=5$}
 \end{array}\right.
 $$
 is a $C_8$ in $G$,  a contradiction (see Figure 11).
 \begin{figure}\label{fig-11}
\psfrag{A}{$x_1$}
\psfrag{B}{$x_2$}
\psfrag{C}{$x_3$}
\psfrag{D}{$x_4$}
\psfrag{E}{$x_5$}
\psfrag{F}{$y_1$}
\psfrag{G}{$z_1$}
\psfrag{H}{$u_1$}
\psfrag{J}{$y_4$}
\psfrag{K}{$z_4$}
\psfrag{L}{$u_4$}
\psfrag{M}{$v_4$}
\psfrag{N}{$k=4$}
\psfrag{U}{$k=5$}
\psfrag{V}{$y_5$}
\psfrag{W}{$z_5$}
\psfrag{X}{$u_5$}
\psfrag{Y}{$v_5$}
      \begin{center}
     \includegraphics[width=130mm]{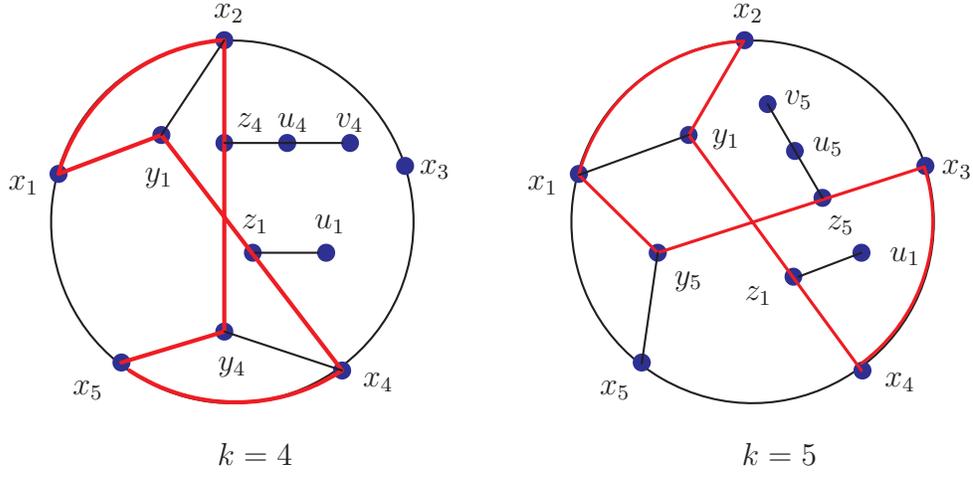}\\
   \caption{The $C_8$ in case $I_C=\{x_1,x_k\}$, $k=4,5$.}
   \end{center}
\end{figure}

\item  $I_C=\{x_1\}$.
Let $x_1y_1z_1u_1v_1$ be a near-good $(C,x_1)$-path satisfying Claim \ref{cla-3-9}. Then, $y_1z_1u_1v_1$ is an induced path of $G-C$ such that
      \begin{equation}\label{eqn-3-47}
      E_G(\{y_1,z_1,u_1,v_1\}, V(C))=\{y_1x_1,y_1x_2,z_1x_4\}.
      \end{equation}
Let $k$ be an integer with $k\in\{3,5\}$. Then, $x_k\notin I_C\cup I_C^+$.
By applying Claim \ref{cla-3-8} with $i=k$, we obtain a good $(C,x_k)$-path, say $x_ky_kz_ku_k$, of length $3$ in $G$. It follows from Definitions \ref{def-2-2} that $y_kz_ku_k$ is an induced path of $G-C$ such that
\begin{equation}\label{eqn-3-48}
E_G(\{y_k,z_k,u_k\}, V(C))= \{y_kx_k\},~k=3,5.
\end{equation}
Note that $y_1\in A_1\cap A_2$ and $y_k\in A_{1+2}\cup A_{1+4}$. By applying Claim \ref{cla-3-5} with $P:=y_1z_1u_1v_1$ and $Q:=y_kz_ku_k$, we can derive that
\begin{equation}\label{eqn-3-49}
\{y_1,z_1,u_1,v_1\}\cap\{y_k,z_k,u_k\}=\emptyset,~k=3,5
\end{equation}
and
\begin{equation}\label{eqn-3-50}
E_G(\{y_1,z_1,u_1,v_1\}, \{y_k,z_k,u_k\})\subseteq \{v_1u_k\},~k=3,5.
\end{equation}
It follows from (\ref{eqn-3-48}) that $y_3\neq y_5$.
Set $H_{35}:=G[\{y_3,z_3,u_3\}\cup\{y_5,z_5,u_5\}]$. If $H_{35}$ is not connected, then by (\ref{eqn-3-48}), we can derive  that
$$
u_5z_5y_5x_5x_1x_2x_3y_3z_3u_3
$$
is an induced $P_{10}$ in $G$, a contradiction (see Figure 12).
\begin{figure}\label{fig-12}
\psfrag{A}{$x_1$}
\psfrag{B}{$x_2$}
\psfrag{C}{$x_3$}
\psfrag{D}{$x_4$}
\psfrag{E}{$x_5$}
\psfrag{F}{$y_3$}
\psfrag{G}{$z_3$}
\psfrag{H}{$u_3$}
\psfrag{I}{$y_5$}
\psfrag{J}{$z_5$}
\psfrag{K}{$u_5$}
\psfrag{L}{$y_1$}
\psfrag{M}{$z_1$}
\psfrag{N}{$u_1$}
\psfrag{O}{$v_1$}
\psfrag{P}{$H_{35}$}
      \begin{center}
     \includegraphics[width=100mm]{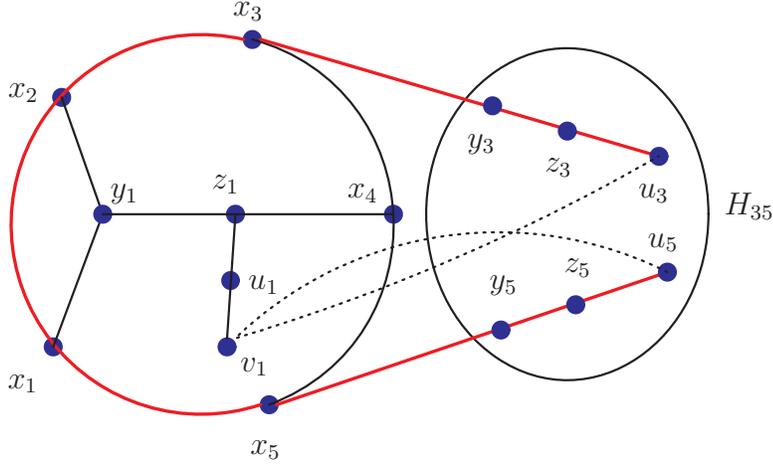}\\
   \caption{Induced $P_{10}$ when $H_{35}$ is not connected}
   \end{center}
\end{figure}
Hence,
\begin{equation}\label{eqn-3-51}
\makebox{$H_{35}$ is a connected subgraph of $G-C$.}
\end{equation}

Let $Q$ be a shortest $(y_3,y_5)$-path in $H$. Then, $\ell(Q)\geq 1$.
If $\ell(Q)=1$, then $y_3y_5\in E(G)$  and $x_3y_3y_5x_5x_1y_1z_1x_4x_3$ is a $C_8$ in $G$, a contradiction. Hence, $\ell(Q)\geq 2$.
Set
$$
R:=\left\{\begin{array}{ll}
   x_3y_3\overrightarrow{Q}y_5x_5x_1y_1z_1u_1v_1, &\makebox{if $u_3,u_5\notin V(Q)$}\cr
   x_3y_3\overrightarrow{Q}y_5x_5x_1y_1z_1u_1, &\makebox{otherwise.}
 \end{array}\right.
$$
By (\ref{eqn-3-47})-(\ref{eqn-3-50}), we see that $R$ is an induced path of $G$. As $G$ is $P_{10}$-free, $\ell(R)\leq 8$. This together with  $\ell(Q)\geq 2$  implies that $\ell(Q)=2$ and $V(Q)\cap (\{u_3\}\cup\{u_5\})\neq\emptyset$. Hence, $Q=y_3u_3y_5$ or $Q=y_3u_5y_5$, contrary to the fact that both $y_3z_3u_3$ and $y_5z_5u_5$ are induced paths of $G-C$.
\end{itemize}
In each case, we get a contradiction. Hence, Claim \ref{cla-3-10} is true.
\qed

\vspace{3mm}
It follows from Claim \ref{cla-3-10} that $I_C=\emptyset$. Let $i$ be an integer with $i\in [1,5]$. Then, $x_i\notin I_C\cup I_C^+$.  By Claim \ref{cla-3-8} and Definition \ref{def-2-2},  we see that $G-C$ contains  an induced path $y_iz_iu_i$  such that
\begin{equation}\label{eqn-3-52}
E_G(\{y_i,z_i,u_i\}, V(C))= \{y_ix_i\}, ~i\in [1,5].
\end{equation}
It follows that $N_C(z_i),N_C(y_1),N_C(y_2),N_C(y_3),N_C(y_4),N_C(y_5)$ are distinct subsets of $V(C)$,
and hence $z_i, y_1,y_2,y_3,y_4,y_5$ are distinct vertices of $G-C$. Set
$$
H_{i}:=G[\{x_i,y_i,z_i\}\cup\{x_{i+2},y_{i+2}, z_{i+2}\}],~i\in [1,5].
$$
By an analogy similar to that in the proof of (\ref{eqn-3-51}), we can derive that
$H_i$ is a connected subgraph of $G-C$. Let $Q_i$ be a shortest $(y_i,y_{i+2})$-path in $H_i$. Then, $\ell(Q_i)\geq 1$. If $\ell(Q_i)= 1$, then $y_iy_{i+2}\in E(G)$. By (\ref{eqn-3-52}), we see that $G[V(C)\cup\{y_i,y_{i+2}\}]\cong\theta(2,3,3)$. This together with Lemma \ref{lem-2-6} implies that $G$ admits a $C_4$ or $C_8$, a contradiction. Hence, $\ell(Q_i)\geq 2$. By applying Claim \ref{cla-3-4} with $(P, j,u,v):=(Q_i,i+2,y_i,y_{i+2})$, we see that $\ell(Q_i)\neq 3,4$.
Hence,
\begin{equation}\label{eqn-3-53}
\ell(Q_i)\in\{2,5\}, ~i\in [1,5].
\end{equation}
We consider the following two cases.

{\bf Case 1.} ~$\ell(Q_i)=5$ holds for some $i\in [1,5]$.

By permuting the vertices of $C$ (if necessary), we may assume that  $\ell(Q_1)=5$. Then,
\begin{equation}\label{eqn-3-54}
E_G(\{y_1,z_1,u_1\},\{y_3,z_3,u_3\})=\{u_1u_3\}.
\end{equation}
This together with (\ref{eqn-3-52}) implies that $y_1,z_1,u_1,y_3,z_3,u_3,y_4$ are distinct vertices of $G-C$. If  $E_G(\{y_4\},  \{y_1,z_1,u_1,y_3,z_3,u_3\})=\emptyset$, then by (\ref{eqn-3-52}) and (\ref{eqn-3-54}), we can derive that $y_4x_4x_3y_3z_3u_3u_1z_1y_1x_1$ is an induced $P_{10}$ in $G$, a contradiction. Hence,
\begin{equation}\label{eqn-3-55}
E_G(\{y_4\},  \{y_1,z_1,u_1,y_3,z_3,u_3\})\neq\emptyset.
\end{equation}
On the other hand, by applying (\ref{eqn-3-53}) with $i=4$, we have $\ell(Q_4)\in\{2,5\}$, and hence $y_4y_1\notin E(G)$. This together with
(\ref{eqn-3-55}) implies that
$y_4w\in E(G)$ holds for some vertex $w\in \{z_1,u_1,y_3,z_3,u_3\}$. Set
$$
C^*:=\left\{\begin{array}{ll}
y_4z_1u_1u_3z_3y_3x_3x_4y_4,&\makebox{if $w=z_{1}$}\cr
y_4u_1z_1y_1x_1x_2x_3x_4y_4,&\makebox{if $w=u_{1}$}\cr
y_4y_3x_3x_4y_4,&\makebox{if $w=y_3$}\cr
y_4z_3y_3x_3x_2x_1x_5x_4y_4,&\makebox{if $w=z_3$}\cr
y_4u_3u_1z_1y_1x_1x_5x_4y_4,&\makebox{if $w=u_3$}
\end{array}\right.
$$
Then, $C^*$ is a $C_4$ or $C_8$ in $G$, a contradiction.

{\bf Case 2.} ~$\ell(Q_i)=2$ for all $i\in [1,5]$.
\begin{figure}\label{fig-13}
\psfrag{a}{$x_1$}
\psfrag{b}{$x_2$}
\psfrag{c}{$x_3$}
\psfrag{d}{$x_4$}
\psfrag{e}{$x_5$}
\psfrag{f}{$y_1$}
\psfrag{g}{$y_3$}
\psfrag{h}{$y_5$}
\psfrag{i}{$y_2$}
\psfrag{j}{$y_4$}
\psfrag{k}{$v_1$}
\psfrag{l}{$v_3$}
\psfrag{m}{$v_5$}
\psfrag{n}{$v_2$}
\psfrag{o}{$v_4$}
      \begin{center}
     \includegraphics[width=60mm]{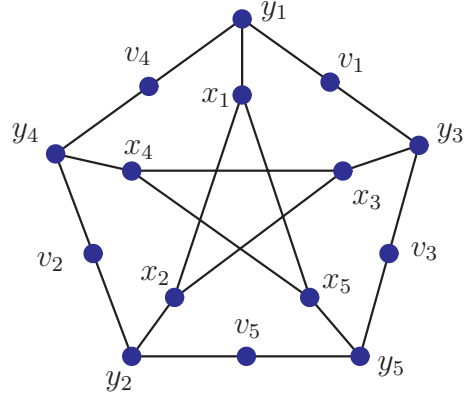}\\
   \caption{The subgraph $H$ of $G$.}
   \end{center}
\end{figure}

For $i\in [1,5]$, let $Q_i:=y_iv_iy_{i+2}$. Then, $v_i\in \{z_i,z_{i+2}\}$. This together with (\ref{eqn-3-52}) implies that
\begin{equation}\label{eqn-3-56}
E_G(\{y_i,v_i\}, V(C))= \{y_ix_i\}, ~i\in [1,5].
\end{equation}
It follows from (\ref{eqn-3-56}) that for each $i\in [1,5]$
\begin{equation}\label{eqn-3-57}
y_1,y_2,y_3,y_4,y_5, v_i ~\makebox{are distinct vertices of $G-C$}.
\end{equation}
We claim that for all $i\in [1,5]$
\begin{equation}\label{eqn-3-58}
N(v_i)\cap \{y_1,y_2,y_3,y_4,y_5\}=\{y_i,y_{i+2}\},
\end{equation}
where $y_6:=y_1$ and $y_7:=y_2$. By way of contradiction, assume that (\ref{eqn-3-58}) is false. Then, $v_iy_j\in E(G)$ holds for some $j\in\{i+1,i+3,i+4\}$. Set $k:=j$ if $j\in\{i+1,i+4\}$, and $k:=j-1$ if $j=i+3$. Then, $x_ky_kv_iy_{k+1}x_{k+1}$ is a $C$-path in $G$.  This together with (\ref{eqn-3-57}) implies that $x_ky_kv_iy_{k+1}x_{k+1}\overrightarrow{C}x_k$
is a $C_8$ in $G$, a contradiction. Hence, (\ref{eqn-3-58}) is true.

It follows from (\ref{eqn-3-57}) and (\ref{eqn-3-58}) that $y_1,y_2,y_3,y_4,y_5,v_1,v_2,v_3,v_4,v_5$ are distinct vertices of $G-C$, and hence the graph $H$, depicted in Figure 13, is a subgraph of $G$. We claim that
\begin{equation}\label{eqn-3-59}
G[\{x_1,x_2,x_3,x_4,x_5,y_1,y_2,y_3,y_4,y_5,v_1,v_2,v_3,v_4,v_5\}]\cong H.
\end{equation}
By way of contradiction, assume that (\ref{eqn-3-59}) is false. Then, by (\ref{eqn-3-56}) and (\ref{eqn-3-58}), we can derive that one of $\{y_1, y_2, y_3, y_4, y_5\}$ and $\{v_1, v_2, v_3, v_4, v_5\}$ is not  independent  in $G$.
And hence there exists a pair of integers $(i,j)$ with $i\in [1,5]$ and $j\in\{i+1,i+2\}$ such that  $\{y_iy_j,v_iv_j\}\cap E(G)\neq\emptyset$. By symmetry, we may assume that $i=1$ and $j\in\{2,3\}$. Then, $G$ admits an edge $e^*$ in $\{y_1y_2,y_1y_3,v_1v_2,v_1v_3\}$. Define
$$
D^*:=\left\{\begin{array}{ll}
y_1y_2x_2x_1y_1,&\makebox{if $e^*=y_1y_2$}\cr
y_1y_3x_3x_2y_2v_2y_4v_4y_1,&\makebox{if $e^*=y_1y_3$}\cr
v_1v_2y_2v_5y_5x_5x_1y_1v_1,&\makebox{if $e^*=v_1v_2$}\cr
v_1v_3y_3x_3x_4x_5x_1y_1v_1,&\makebox{if $e^*=v_1v_3.$}
\end{array}\right.
$$
Then, $D^*$ is a $C_4$ or $C_8$ in $G$, a contradiction. Hence, (\ref{eqn-3-59}) is true.

It follows from (\ref{eqn-3-59})  that $x_1y_1v_1y_3v_3y_5v_5y_2v_2y_4$ is an induced $P_{10}$ in $G$, a contradiction. This completes the proof of Theorem \ref{thm1.1}.
\qed
\\










\end{document}